\documentclass[12pt]{article}

\usepackage[latin1]{inputenc}
\usepackage[english]{babel}
\usepackage{amsfonts}
\usepackage{euscript}
\usepackage{geometry, graphicx, pstricks, xypic}

\newtheorem{conj}{Conjecture}
\newtheorem{theo}[conj]{Theorem}
\newtheorem{prop}[conj]{Proposition}
\newtheorem{coro}[conj]{Corollary}
\newtheorem{lemm}[conj]{Lemma}

\newtheorem{defi}[conj]{Definition}
\newtheorem{nota}[conj]{Notation}
\newtheorem{exem}[conj]{Example}
\newtheorem{rema}[conj]{Remark}

\newcommand{\bx}{{\bf x}}
\newcommand{\by}{{\bf y}}

\newcommand{\Spec}{{\rm Spec}\ }
\newenvironment{proo}{{\flushleft \it Proof.}}{\hfill $\square$ \vspace{2mm}}

\renewcommand{\choose}[2]{
\left (
\begin{array}{c}
\hspace{-2mm}  #2  \hspace{-2mm} \\
\hspace{-2mm}  #1  \hspace{-2mm}
\end{array}
\right )
}

\def \cI {{\cal I}}

\def \mC {{\mathfrak C}}

\newcommand{\vacuum}{\phi}
\newcommand{\sn}{{S^{[n]}}}
\newcommand{\inc}{S^{[n,n+1]}}
\newcommand{\cinc}{(\p^2)^{[n,n+1]}}
\newcommand{\inco}{S^{[n,n+1]}_0}
\newcommand{\cinco}{(\p^2)^{[n,n+1]}_0}
\newcommand{\ch}[1]{A_T^*(#1)}
\newcommand{\ck}[1]{A_K^*(#1)}
\newcommand{\Ia}{{\cal I}_{1i}}
\newcommand{\Ib}{{\cal I}_{i1}}

\newcommand{\fix}{\mathrm{fix}}
\newcommand{\es}{\mathrm{es}}
\newcommand{\nak}{\mathrm{nak}}
\newcommand{\Cok}{\mathrm{Coker}}
\newcommand{\Ker}{\mathrm{Ker}}
\newcommand{\Tan}{\mathrm{Tan}}
\newcommand{\scal}[1]{\langle #1 \rangle}
\newcommand{\ES}{Ellingsrud-Stromme }
\newcommand{\SD}[1]{S^{[#1]}_\Delta}
\newcommand{\SDL}[2]{S^{[#1]}_{\Delta,#2}}

\newcommand{\N}{\mathbb{N}}

\newcommand{\lt}[1]{

\noindent
\textbf{LT:}
#1

\noindent
}
\newcommand{\comm}[1]{}

\newcommand{\im}{\mathrm{Im}\ }
\renewcommand{\Im}{\mathrm{Im}}

\newcommand{\coim}{\mathrm{coIm}}
\newcommand{\supp}{\mathrm{supp}}
\let \fl=\flushleft

\let \beg=\begin
\let \mc=\mathcal

\newcommand{\fb}{\ensuremath{\downarrow}}

\newcommand{\fd}{\ensuremath{\rightarrow}}

\newcommand{\pla}{\A^2}

\newcommand{\x}{\ensuremath{\times}}
\newcommand{\pn}{{\cal P}_n}

\newcommand{\p}{{\mathbb P}}
\newcommand{\A}{{\mathbb A}}
\newcommand{\Q}{{\mathbb Q}}
\newcommand{\Z}{{\mathbb Z}}

\newcommand{\lpara}{\vskip 3mm}

\begin{document}

\title{On the equivariant cohomology of Hilbert schemes of points in the plane}
\author{Pierre-Emmanuel Chaput, Laurent Evain}
\maketitle

\section*{Abstract}
\label{sec:abstract}
\comm{
Let $\pla$ be the affine plane
regarded as a toric variety with an action of the torus $T$.
We study the equivariant Chow ring $A_{T}^*(Hilb^n(\pla))$
of the punctual Hilbert scheme $Hilb^n(\pla)$, and its extension
$A_{T,Q}^*(Hilb^n(\pla))=A_{T}^*(Hilb^n(\pla))\otimes_{A_T^*(point)}Q$,
with $Q$ the fraction field of $A_T^*(point)$.
We introduce new
operators acting on the direct sum
$\oplus_{n\in \mathbb{N}}A_{T,Q}^*(Hilb^n(\pla))$. We give formulas
to compute both the classical operators (creation operators $q_i$, intersection
with the boundary...), these new operators and some commutation
relations between them. We give a formula for the
class in $A_T^*(Hilb^n(\pla))$ of the small diagonal in $Hilb^n(\pla)$
in terms of the (equivariant) Chern class of the tautological bundle.
We compute base change formulas  in  $A_{T,Q}^*(Hilb^n(\pla))$
between the natural basis
introduced by Nakajima, Ellingsrud and Str\o mme, and the
classical basis associated with the fixed points.
Finally, we compute the equivariant Chow ring  $A_T^*(Hilb^n(\pla))$,
generalizing the work of Vasserot.
}

Let $S=\pla$ be the affine plane
regarded as a toric variety with an action of the 2-dimensional torus $T$.
We study the equivariant Chow ring $A_{K}^*(\sn)$
of the punctual Hilbert scheme $\sn$ with equivariant coefficients inverted.
We compute base change formulas in $A_{K}^*(\sn)$
between the natural bases
introduced by Nakajima, Ellingsrud and Str\o mme, and the
classical basis associated with the fixed points. We compute the
equivariant commutation relations between creation/annihilation operators.
We express the class of the small diagonal in $\sn$ in terms of the
equivariant Chern classes of the tautological bundle. We prove that the 
nested Hilbert scheme $\inco$ parametrizing nested punctual subschemes of 
degree $n$ and $n+1$ is irreducible.

\section*{Introduction}

If $S$ is a quasi-projective smooth surface, let $S^{[n]}$ be the
Hilbert scheme parameterizing the zero dimensional subschemes of
degree $n$ in $S$.
Following Nakajima and Grojnowski, a first tool to
study the Chow ring $A^*(S^{[n]},\Q)$
is to consider the direct sum $\oplus_{n\in \N} A^*(S^{[n]},\Q)$
and
operators acting linearly on this direct sum. Then, a lot of structure
and information lies in the commutation
relations of the various operators. In the case $S=\pla$, this
approach yields a basis of $A^*(S^{[n]},\Q)$ that we call
Nakajima's basis and a description of the
ring structure on it \cite{nakajima,lehn}.

When $S=\pla$, another approach is the use of the equivariant
Chow rings. The 2-dimensional torus $T$ acts on $S^{[n]}$.
The equivariant Chow ring with respect to the action of the full
torus $T$ has been computed in \cite{laurent} in the case
$S=\mathbb{P}^2$, but this is a purely equivariant approach independent of
Nakajima's framework. Similarly Bialynicki-Birula's theorem \cite{BB}
yields a basis of the classical and equivariant Chow rings which has been
studied in \cite{ES} and which we call Ellingsrud-Str\o mme's basis.

There are equivariant
analogues of the operators introduced by Nakajima et al
which act on the equivariant Chow ring. Following Vasserot
\cite{vasserot}, it is natural
to compute these equivariant operators. In his paper, Vasserot does
not consider the full action of the torus $T$, but the action of a
non-generic one-dimensional subtorus $T' \subset T$. He computes several
operators in $T'$-equivariant Chow rings and their commutators. As a
consequence, he obtains a description of the $T'$-equivariant and of
the classical Chow ring of the Hilbert scheme of $\pla$.

On the other hand, Schiffmann and Vasserot study an algebra of operators
acting on the equivariant $K$-theory of Hilbert schemes \cite{SV}. Since
the correspondences defining $q_i$ for $i>1$ are singular, they do not
define operators on the $K$-theory, this is the reason why the authors
only consider the operator algebra generated by $q_1,q_{-1}$ and multiplication
by some tautological bundles.

\lpara

In this work we consider the action of the operators $q_i$ for all $i$
on the $T$-equivariant Chow rings.
The apparent difficulty coming from the non projectivity of $\A^2$
is not severe : we have all the standard constructions and properties
of intersection theory that we need (pushforward, correspondences,
composition of correspondences...) provided that we work in the
tensored equivariant Chow ring
$A_{T}^*(\sn) \otimes_{A_T^*(pt)} K$,
instead of $A_{T}^*(\sn)$, where $K$ is the fraction field
of $A_T^*(pt)$, the equivariant Chow ring of a point
(Section \ref{sec:pushforward-with-non}).

However this construction also has its drawback: the pushforward of
a contractant non proper morphism does not need to vanish (see Lemma
\ref{lemm:bott} for an example) and some key arguments of the classical
situation are not valid in our equivariant context.
Let $S^{[n]}_0$ denote the set of subschemes $z_n$ of length $n$ supported as the
origin, and let $S^{[p,q]}_0$ denote the similar set of couples
of nested subschemes $(z_p \subset z_q)$. 
When one wants to compute the composition of two correspondences,
the ubiquitous local situation that one has to understand are the classes
$\pi_{*} [C]$, where $C$ is some subvariety
in $S^{[p,q]}$ and $\pi$ is the projection to $S^{[p]}$ or to $S^{[q]}$. 
The geometry  is under control 
when both $z_p$ and $z_q$ are curvilinear for the generic pair $(z_p,z_q)\in C$.
In the other cases, the restriction of $\pi$
to $C$ is contractant and therefore $\pi_* [C]$ is zero in the classical
Chow ring \cite{nakajima,lehn}. However $\pi_* [C]$  need not vanish in the equivariant
Chow ring.
Our remedy is to prove that $S^{[n,n+1]}_0$ is irreducible and well understood
(Theorem \ref{theo:inco}, Section \ref{sec:tangent-space}) (whereas
$S^{[p,q]}_0$ is not in general, Proposition \ref{pro:s24}).
Then follows our construction to compute the commutators:
we use algebraic arguments to reduce to the
case when one of the operators adds only one point.

\lpara

In Section \ref{sec:classical}, we consider the classical operators
acting on $A_K:=\oplus_{n\in \N}  A_{T}^*(\sn) \otimes_{A_T^*(pt)} K $, 
namely the creation/destruction operators $q_i$,
the boundary operator $\partial$, and an auxiliary operator $\rho$.
All these operators are defined by a correspondence. Provided that
the correspondence is smooth, the computation is easily done with
the Bott formula. This is the strategy to compute $q_1$ and
$q_{-1}$ in the fixed point basis (Proposition \ref{pro:q1} and \ref{pro:q-1}). 
All the other correspondences are singular at some
points and a turnaround is needed to compute the corresponding operators.

Computing restriction to fixed points, we prove the formula
$\partial = -2 c_{1}({\cal O}_X[n])$, where
${\cal O}_X[n]$ denotes the tautological bundle.
Following Lehn and Schiffmann-Vasserot's ideas, we
consider various commutators
starting with $q_1,q_{-1}$ and $\partial$. We end up with recursion formulas
for the $q_i$'s, $|i|>1$ (Theorem \ref{theo:qi}).
In particular this yields 
base change formulas between the fixed point basis and Nakajima's
basis (Example \ref{ex:nakToFix}).

To compute the commutation relations
between the $q_i$'s (Theorem \ref{theo:qi-qj}), using once again
the same general
idea as in \cite{SV}, we use algebraic computations to reduce to the
case of operators of conformal degree one. The algebraic 
reduction leads to Proposition
\ref{pro:rho-rhod},
a formula apparently new even in the non equivariant context. We use
geometric arguments to get rid of the excess intersection components 
which appear in this case. 
\comm{
To compute $\partial$, we basically need to compute the restriction
$\Delta_p$ of the diagonal $\Delta$ at a (possibly singular) point
$p$. We introduce a smooth variety $F$ such that the intersection
$\Delta \cap F$ is a singular union of smooth varieties $S_i$. In
particular the restriction $(\Delta \cap F)_p=\sum (S_i)_p$ is computed
easily and we conclude with the formula $\Delta_p= (\Delta \cap
F)_p/F_p$.
}

The class $\delta_n \in A^{n-1}(\sn)$ of the small
diagonal $\Delta_n\subset \sn$ parameterizing the subschemes
supported on a single point has an expression in terms of the equivariant Chern classes of the
tautological bundle: $\delta_n=(-1)^{n-1}n c_{n-1}({\cal O}_X[n])$.
The originial proof by Lehn \cite[Theorem 4.6]{lehn}
remains true in our context. We give a new proof which relies on 
an algebraic expression for the operator $q_n$ 
(Theorem \ref{theo:calculQn}).
\lpara

Finally, we give an application of our equivariant computations 
at the level of classical Chow rings. We investigate the base change between
Nakajima's basis and Ellingsrud
and Str\o mme's basis. Our theorem asserts that these bases are
equal up to sign (and a normalisation constant with our conventions) 
in classical Chow rings (Theorem \ref{theo:nak2es}).  Our strategy of proof is
to interpret the Bialynicki-Birula cells in terms of operators:
we introduce new creation
operators $q_{i,X}$ such that the basis introduced  by 
Ellingsrud and Str\o mme is obtained
applying these operators on the vacuum. We express the
$q_{i,X}$ in terms of the creation operators $q_i$ and we get a base
change formula in the equivariant Chow ring. Projecting this
relation in the usual Chow ring gives the asserted formula.

\lpara

\textit{Acknowledgements: We thank the Institut Henri Poincar\'e in
  Paris and the Mathematisches Forschungsinstitut Oberwolfach
where part of this reaserch took place. 
We thank the developers of the Sage project for their software.}

\tableofcontents

\section{Pushforward with non proper morphisms}
\label{sec:pushforward-with-non}

We work over an algebraically closed field $k$ of any characteristic.
Let $T$ be a 2-dimensional torus. The $T$-equivariant Chow ring
$A_T^*(pt)$ of a point is isomorphic to a polynomial
ring in two variables $U,V$.
We denote by $K=\Q(U,V)$ the field of fractions of $A^*_T(pt)$.
Moreover, if $X$ is any $T$-variety, we denote by $A^*_K(X)$ the tensor
product $A^*_T(X) \otimes_{\Z[U,V]} K$. We denote the product of two
classes $x,y$ in a Chow
ring indifferently by $x.y$ or by $x\cup y$.

In this section $f:X \to Y$ is an equivariant morphism between smooth
varieties. Moreover, we assume that $X$ and $Y$ are {\em filtrable}, in the
sense of Definition 3.2 in \cite{brion}.

When $f:X \fd Y$ is a proper equivariant morphism,
there is a well defined pushforward
$f_*:A_T^*(X) \fd A_T^*(Y)$.  Since we shall work with the affine plane,
we are in a non projective setting and we have to deal with non
proper morphisms.

The goal of this section is to explain that a good notion of
pushforward $f^K_{*}$ exists, when $f$ is a non proper morphism,
provided that the restriction to fixed points $f^T:X^T\fd Y^T$ is
proper.  This notion is applied to define correspondences. We show
that these correspondences defined using non proper morphisms
satisfy formal properties similar to the correspondences defined in
the usual setting when $f$ is proper.

\begin{defi}
\label{def:pdt-tensoriel}
  If $f$ is as above, $f^T:X^T\fd Y^T$ will denote the
  restriction of $f$ to $T$-fixed points. The morphisms
  $f^K_* : A_K^*(X) \fd A_K^{*+\dim Y-\dim X}(Y)$
  (when $f$ is proper)
  and $f_K^*:A_K^*(Y) \fd A_K^*(X)$
  are derived from the standard morphisms
  $f^T_*: A_T^*(X) \fd A_T^{*+\dim Y-\dim X}(Y)$ and
  $f^*_T:A_T^*(Y) \fd A_T^*(X) $
  after tensorisation over $\Z[U,V]$ by $K$.
\end{defi}

Let $f:X \to Y$ be any $T$-equivariant morphism, and
consider the following commutative diagram:
\begin{equation}
\label{equ:i}
\begin{array}{cccl}
  X^T     & \stackrel{i}{\hookrightarrow}   & X \\
  \fb f^T &                               & \fb f \\
  Y^T     & \stackrel{j}{\hookrightarrow} & Y & .
\end{array}
\end{equation}
Since $i^K_*$ is an isomorphism by \cite[Theorem 1]{EG},
the following definition is meaningful:
\begin{defi}
\label{def:propre}
  If $f^T$ is proper, define
  $$f^{T,K}_* = j^K_* (f^T)^K_* {(i_*^K)}^{-1} : A^*_K(X)
  \to A^{*+\dim Y-\dim X}_K(Y) \ .$$
\end{defi}
\noindent
If $f$ is proper, then $f^{T,K}_* = f^K_*$, by the functoriality of
the pushforward in the proper case. Since there is therefore no possible
confusion, we will denote $f^{T,K}_*$ simply by $f^K_*$.

\begin{theo}
\label{theo:fk*}
  The morphism $f^K_*$ satisfies the following properties:
  \begin{enumerate}
  \item Functoriality: if we have $T$-equivariant morphisms
    $X \stackrel{f}{\to} Y \stackrel{g}{\to} Z$ such that $f^T$ and $g^T$
    are proper, then
    $(g\circ f)^K_* = g^K_* \circ f^K_*$.
  \item Projection formula: assume here that $X$ and $Y$ are smooth, so that
    $A^*_K(X)$ and $A^*_K(Y)$ are rings.
    For any $\alpha \in A^*_K(X)$ and
    $\beta \in A^*_K(Y)$, we have the equality
    $f^K_* (\alpha) \cdot \beta = f^K_* ( \alpha \cdot f_K^*(\beta) )$.
  \item We have the equality $\ g^K_* f^*_K = l^*_K h^K_*\ $
    if $f,g,h,l$ are as in the following diagram:
    $$
    \begin{array}{ccc}
      X\x Y \x Z & \stackrel{g}{\fd} & Y \x Z \\
      \fb f      &                   & l\fb \\
      X\x Y      & \stackrel{h}{\fd} & Y
    \end{array}
    $$
\end{enumerate}
\end{theo}
\begin{proo}
{\fl \textbf{Functoriality.}} Consider the following diagram:
$$
  \begin{array}{cccl}
  X^T     & \stackrel{i}{\hookrightarrow}   & X \\
  \fb f^T &                               & \fb f \\
  Y^T     & \stackrel{j}{\hookrightarrow} & Y \\
  \fb g^T &                               & \fb g \\
  Z^T     & \stackrel{k}{\hookrightarrow} & Z & .
\end{array}
$$
Then we have
\begin{eqnarray*}
(g\circ f)^K_*&=&k^K_*((g\circ f)^T)^K_*(i^K_*)^{-1}\\
&=& k^K_*(g^T)^K_* (f^T)^K_* (i^K_*)^{-1}\\
&=& k^K_* (g^T)^K_* (j^K_*)^{-1}\ \circ \ j^K_* (f^T)^K_* (i^K_*)^{-1}\\
&=& g^K_* \circ f^K_* \ \ .
\end{eqnarray*}
\\
\textbf{Projection formula}. Let $\alpha\in A^{*}_K(X)$ and
$\beta\in A^{*}_K(Y)$. By commutativity of the diagram (\ref{equ:i}), we have:
\begin{eqnarray*}
\alpha \cdot f^*_K \beta & = & \alpha \cdot {(i_K^{*})}^{-1}
(f^T)_K^* j^*_K \beta\\
&=&  i^K_* {(i^K_*)}^{-1} \alpha \cdot {(i_K^*)}^{-1} (f^T)_K^* j^*_K \beta\\
&=&  i^K_* ( {(i^K_*)}^{-1} \alpha \cdot
i_K^* {(i_K^*)}^{-1} (f^T)_K^* j^*_K \beta),\\
\end{eqnarray*}
where the last equality is due to the usual projection formula for $i^K_*$.
Applying $(i^K_*)^{-1}$ to both sides of the equality we get:
$(i^K_*)^{-1}(\alpha \cdot f^*_K\beta)=
(i^K_*)^{-1} \alpha \cdot (f^T)_K^* j^*_K \beta.$
Therefore,
\begin{eqnarray*}
f^K_* \alpha \cdot \beta & = & j^K_* (f^T)^K_*
(i^K_*)^{-1} \alpha \cdot \beta\\
& = & j^K_* (f^T)^K_* \, ( \, (i^K_*)^{-1} \alpha
\cdot (f^T)^*_K j^*_K \beta \, )\\
& = & j^K_* (f^T)^K_* \, ( \, (i^K_*)^{-1} (\alpha \cdot f^*_K \beta) \, )\\
& = & f^K_* (\alpha \cdot f^*_K\beta )
\end{eqnarray*}
(the second equality follows from the projection formula for $j \circ f^T$.)
\\
\textbf{Property 3.} Consider the following diagram:
$$
    \begin{array}{cccccccl}
      X^T\x Y^T \x Z^T & \stackrel{b}{\hookrightarrow}& X\x Y \x Z &
      \stackrel{g}{\fd} & Y \x Z & \stackrel{d}{\hookleftarrow} &Y^T \x
      Z^T \\
      {f^T}\fb &&\fb f      &                   & l\fb && \fb l^T\\
      X^T \x Y^T & \stackrel {a}{\hookrightarrow} & X\x Y      &
      \stackrel{h}{\fd} & Y & \stackrel{c}{\hookleftarrow} & Y^T & .
    \end{array}
$$
The property that we want to prove holds when the pushforward maps
are taken with proper morphisms, so that for example $b^K_*
(f^T)^*_K = f_K^* a^K_*$. The second following equality is a
consequence of this remark, the third and the fourth equalities are similar:

\begin{eqnarray*}
g^K_*f^*_K & = & d^K_* (g^T)^K_* {(b^K_*)}^{-1} f^*_K \\
& = & d^K_* (g^T)^K_* (f^T)_K^* {(a^K_*)}^{-1} \\
& = & d^K_* (l^T)_K^* (h^T)^K_* {(a^K_*)}^{-1} \\
& = & l^*_K c^K_* (h^T)^K_* {(a^K_*)}^{-1} \\
& = & l^*_K h^K_* \ .
\end{eqnarray*}
\end{proo}

In practice, $f^K_*$ can be computed by a ``Bott formula'', as in the
proper case. Assume that $X$ is smooth.
Since $X^T$ is smooth, $A_T^*(X^T)=\oplus A_T^*(X_i)$ where
the sum runs through the irreducible components $X_i$ of $X^T$.
We denote by $c_{top}(N_{X^T,X})$ the operator which acts
on $A_T^*(X^T)$
 through
multiplication by the equivariant Chern class $c_{d_i}$ of the normal bundle
$N_{X_i,X}$ on the component $A_T^*(X_i)$, where $d_i$ is the codimension
of $X_i$ in $X$. Similarly, there
is a class $c_{top}(N_{Y^T,Y})$. In $A^*_K(X)$ (or $A^*_K(Y)$), the Chern class
$c_{d_i}$ is equal to the sum of an invertible element and a nilpotent element,
according to the proof of \cite[Proposition 3.2(i)]{brion}.
Therefore it is invertible and $c_{top}(N_{X^T,X})$ is an invertible operator.

\begin{lemm}
\label{lem:i*}
Assume that $X$ is smooth.
Let $i:X^T \to X$ be the natural inclusion.
The pullback $i_K^*$ is invertible with inverse
$i_*^K \frac{1}{c_{top}(N_{X^T,X})}\ $.
\end{lemm}
\begin{proo}
For $\alpha \in A^*_T(X^T)$, the self-intersection formula reads
$$
i_T^*\, i_*^T\, \alpha = c_{top}(N_{X^T,X})\, \alpha\ .
$$
Thus the lemma follows from the fact that $i_*^K$ is
invertible by the localization theorem \cite[Theorem 1]{EG}. 
\end{proo}

\begin{theo}[Bott Formula]
\label{theo:bott}
  Recall the diagram (\ref{equ:i}) and assume that $X$ and $Y$ are smooth.
  Let $\alpha \in A_T^*(X)$. Then
  $$j^*_Kf^K_* (\alpha) = c_{top}(N_{Y^T,Y})\
  (f^T)^K_* \left (\frac{1}{c_{top}(N_{X^T,X})}
  i_{K}^*(\alpha) \right ).$$ In particular, when
  both $X$ and $Y$ have a finite number of fixed points
  $x_1,\dots,x_n$, $y_1,\dots,y_p$, the formula expresses $f^K_*$ in
  terms of the localization at these fixed points:
  $$(f^K_* \alpha)(y_k)=\sum_{f(x_i)=y_k}\frac{c_{top}(T_{y_k,Y})}{c_{top}
    (T_{x_i,X})} \alpha_{x_i}\ ,$$ where $T_{y_k,Y}$ and $T_{x_i,X}$
  are the tangent $T$-representations.
\end{theo}
\begin{proo}
  By Lemma \ref{lem:i*}, the inverse of $i^K_*$ is
  $\frac{1}{c_{top}(N_{X^T,X})}i^*_K $. Then,
  \begin{eqnarray*}
    j^*_Kf^K_* \alpha & = & j^*_Kj^K_* (f^T)^K_* {(i^K_*)}^{-1} \alpha \\
    & = & c_{top}(N_{Y^T,Y})(f^T)^K_* {(i^K_*)}^{-1} \alpha \\
    & = & c_{top}(N_{Y^T,Y})\ (f^T)^K_* \left ( \frac{1}{c_{top}(N_{X^T,X})}
    i_{K}^*(\alpha) \right ).
  \end{eqnarray*}
\end{proo}

\begin{defi} An equivariant correspondence is a closed $T$-stable
  subvariety  $C\subset X\x Y$ such that $C^T \fd Y^T$ is a proper
  morphism. Let $\pi_X$ and $\pi_Y$ be the projections from $C$ to $X$
  and $Y$ respectively.
  The classes of such varieties $C$ generate a subspace in $A_K^*(X\x Y)$
  and we still call equivariant correspondence a class in this
  subspace.
  An equivariant correspondence $C$
  yields a morphism $f:A_K^*(X) \fd A_K^*(Y)$ defined by
  $f(\alpha)= (\pi_Y)_*^K (\pi_X)_K^*(\alpha))$.
\end{defi}

\begin{prop}
  Assume that $X,Y,Z$ are smooth varieties.
  Let $C\subset X\x Y$ and $D\subset Y\x Z$ be two equivariant
  correspondences, and $f$ and $g$ the associated morphisms. Let
  $\pi_{12},\pi_{13},
  \pi_{23}$ the projections from $X\x Y \x Z$ to $X\x Y,X\x Z$ and
  $Y\x Z$ respectively. Then the
  $(\pi_{13})^K_*((\pi_{12})_K^*[C]\, \cup (\pi_{23})_K^*\, [D])$
  is an equivariant correspondence with associated morphism $g\circ f$.
\end{prop}

\begin{proo}
  When $f:C \fd Y$ and $g:D\fd Z$ are proper, the proof is classical
  and relies on functoriality of the pushforward, on the projection
  formula, and on the third
  property of Theorem \ref{theo:fk*}. Thus the proof is
  mutatis mutandis the same as in the proper case
  if we check that  $E=(\pi_{13})^K_*((\pi_{12})_K^*[C]\, \cup
  (\pi_{23})_K^*\, [D])\in A_K^*(X\times Z)$ is
  well defined and if we can associate to $E$ a morphism $f_E:A_K^*(X)
  \fd A_K^*(Z)$ with the help of the standard formulas of
  intersection theory.

Let $L=(\pi_{12}^{-1}C\cap \pi_{23}^{-1}D)$. Since $D^T \fd Z^T$ is
proper, so is $X^T\x D^T \fd X^T \x Z^T$.
It follows that $L^T=(C^T
\x Z^T)\cap (X^T\x D^T)\fd X^T\x Z^T$ is proper and that
$(\pi_{13})^K_*(L)$ is well defined.

We now prove that $h:L^T \fd Z^T$ is proper.  In particular, there is a
well defined morphism  $(\pi_3^L)_*^K:A_K^*(L) \fd A_K^*(Z)$.
We consider the factorization of $h:L^T\fd Z^T$ as $g\circ f$ with
$f:L^T \fd D^T$ and $g:D^T\fd Z^T$. The morphism $g$ is proper since
$D$ is an equivariant correspondence. The morphism $f:L^T \fd D^T$
is the extension of the proper morphism $C^T \fd Y^T$ by the
morphism $D^T  \fd Y^T$. Then $h$ is proper as a composition of
proper morphisms.

When the intersection $L$ is proper, then $E={(\pi_{13})}^K_*(L)$
and $f_E$ is defined to be  the composite morphism
${(\pi_3^L)}_*^K \circ {(\pi_1^L)}^*_K: A_K^*(X)
  \fd A_K^*(L) \fd A_K^*(Z).$

When the intersection $L$ is not proper, then
${(\pi_{12})}_K^*[C]\, \cup
  {(\pi_{23})}_K^*\, [D]$ is representable by a refined
intersection $I$, ie. by a linear combination of classes $I=\sum
\alpha_i [L_i]$, with $L_i\subset L$  a closed subvariety. In particular,
$L_i^T \fd Z^T$ is proper as a composition of the proper morphisms
$L_i^T\fd L^T \fd Z^T$. Similarly $L_i^T \fd X^T \times Z^T$ is proper
since it factorizes through $L^T$. Then  the pushforward
$E = {(\pi_{13})}^K_{*}( \sum \alpha_i L_i )$  is well
defined. The associated morphism is $f_E=\sum \alpha_i f_{L_i}$ with
$f_{L_i}=  {(\pi_3^{L_i})}_*^K \circ {(\pi_1^{L_i})}^*_K.$
\end{proo}

\begin{defi}
\label{def:produit-scalaire}
Suppose that $\pi:X^T \to \Spec k$ is proper.
Then there is a well-defined pushforward map $\pi^K_* : A_K^*(X) \fd K$. This yields a
$K$-bilinear product on  $A_K^*(X)$ defined by
$\langle \alpha,\beta \rangle_X = \pi^K_*(\alpha \cup \beta)$.
\end{defi}
When $\langle .,. \rangle_X$ and $ \langle .,. \rangle_Y$
are both non degenerate (for instance when $X$ and
$Y$ have a finite number of fixed points), then every map
$f:A_K^*(X) \fd  A_K^*(Y)$ admits a dual map
$f^\vee: A_K^*(Y) \fd  A_K^*(X)$.

\begin{defi}
\label{def:corr-duale} If $C\subset X\x Y$ is a correspondence, the
dual correspondence $C^\vee$ is the correspondence in $Y\x X$ which
is canonically identified with $C$ under the natural isomorphism $X
\x Y \simeq Y\x X$. In particular, if $C^\vee$ is an equivariant
correspondence, it yields a map $A_K^*(Y) \fd  A_K^*(X)$.
\end{defi}

\begin{prop}
\label{pro:duale}
  Assume that $X$ and $Y$ are smooth.
  Let $C\subset X\x Y$ be an equivariant correspondence and
  $f:A_K^*(X) \fd  A_K^*(Y)$ the associated morphism. Suppose that
  $C^\vee\subset Y\x X$ is an equivariant correspondence and that
  $\langle .,. \rangle_X$ and $ \langle .,. \rangle_Y$
  are non degenerate. Then the dual map $f^\vee$ is defined by the dual
  correspondence $C^{\vee}$.
\end{prop}
\begin{proo}
  In the classical setting, the proof relies on the functoriality of
  the pushforward and on the projection formula. Both arguments remain
  valid in our context.
\end{proo}

\noindent
Restriction to fixed points does not commute with the bilinear
product (see Lemma \ref{lem:norme-fix}).
This remark is important when one wants to compute $f^\vee$
on fixed points.

\section{Tangent space to $\inc$}
\label{sec:tangent-space}
We denote by $S$ the affine plane $\A^2$. Let $n \geq 0$ be an
integer, we denote by $S^{[n]}$ the Hilbert scheme parameterizing
length $n$ subschemes of $S$ \cite{grothendieck}
($S^{[0]} = \mathrm{Spec}\ k$). Given $z \in
S^{[n]}$, we denote by $I_z \subset k[X,Y]$ the corresponding
ideal, of codimension $n$. Given $p,q$ integers with
$0 \leq p < q$,
we denote by $S^{[p,q]}$ the ``nested'' Hilbert scheme, namely
the subscheme of $S^{[p]} \times S^{[q]}$ consisting of pairs
$(s,b)$ such that $I_s \supset I_b$. The torus ${(k^*)}^2$ will be
denoted by $T$. It acts on the plane $S$: we use the convention that
an element $(u,v) \in T$ acts on a monomial $X^aY^b \in k[X,Y]$ by
$(u,v) \cdot X^aY^b = (uX)^a(vY)^b$. This induces an action of $T$
on each Hilbert scheme $S^{[n]}$ and $S^{[p,q]}$. We will denote by
$aU+bV$ the weight on $T$ defined by $(aU+bV)(u,v) = u^a v^b$. Given
a monomial $m=X^aY^b$, we denote by $wt(m)=aU+bV$ its weight.
Any character
of $T$ defines naturally an element in $A^1(pt)$, thus our notation is
compatible with the notation $A^*(pt) = \Z[U,V]$ in Section 
\ref{sec:pushforward-with-non}.

{\fl Several} arguments in the present paper rely on a tangent space
argument. In fact, at a $T$-fixed point $z \in S^{[n]}$, the tangent
space $T_z S^{[n]}$ has several combinatorial descriptions. One
\cite{irreductible} is in terms of significant cleft pairs and
another \cite{nakajima} in terms of boxes of the corresponding
staircase. We recall in this section the necessary material to be
comfortable with these two notions. We give two applications. First,
we compute the tangent space at a toric point in $\inc$ as a
representation of $T$. Then Theorem \ref{theo:inco}
proves the irreducibility of $S_0^{[n,n+1]}\subset
S^{[n,n+1]}$ parameterizing the pairs of subschemes $s\subset b$
with the support of $b$ equal to the origin.

\subsection{Tangent space to the Hilbert schemes}

\label{subsection:tangent}

First, observe that a $T$-fixed point $z$ in $\sn$ is defined by an ideal
$$
I_z = \bigoplus_{(a,b) \not \in E} k \cdot X^aY^b\ \ ,
$$
where $E\subset \N^2$ satisfies $(\N^2 \setminus E)+\N^2
\subset (\N^2 \setminus E)$. 
Such (finite) subsets $E\subset \N^2$ will be
called {\em staircases}.
A partition $\lambda = (\lambda_1 \geq \ldots \geq
\lambda_l > 0)$ is by definition a finite sequence of non increasing
positive natural numbers, $l$ is called the length of $\lambda$, and
$|\lambda|=\sum_{i=1}^l \lambda_i$ is the weight of $\lambda$. 
We denote by $\pn$ the set of partitions of
weight $n$. If $E$ is a finite staircase associated with a $T$-fixed
point $z\in \sn$, there exists a unique partition $\lambda$ with
weight $n$ such that 
$(a,b)\in E \Leftrightarrow a+1\leq l, b<\lambda_{a+1}$.

{\fl We begin} by recalling the description given in \cite{irreductible}
of $T_z \sn$ when $z \in \sn$ is a $T$-fixed point.

\begin{itemize}
\item
A monomial $c \in I_z$ is called a {\it cleft} whenever
$X^{-1} \cdot c \not \in I_z$ and $Y^{-1} \cdot c \not \in I_z$.
\item
A Laurent monomial is called {\it positive} (resp., {\it negative})
if it belongs to $Y^{-1}k[X,Y^{-1}]$ (resp., $X^{-1}k[X^{-1},Y]$).
\item
A weight $aU+bV$ with $a \geq 0$ and $b<0$ resp. $a<0$ and $b \geq 0$
will be called {\it positive} resp. {\it negative}.
\item
A {\it cleft pair} is a pair $(c,m)$ such that $c$ is a cleft,
$m$ is a monomial not belonging to $I_z$, and $m/c$ is either positive or
negative (in which case, we say that $(c,m)$ is positive or negative,
respectively).
\end{itemize}

{\fl Now let} $\mc{C}:= \{c_1,\ldots,c_l\}$ denote the set of clefts, which
we order following the convention that $c_{i+1}/c_i$ be positive for
$1 \leq i \leq l-1$.

{\fl For each} positive (resp., negative) cleft pair $(c_i,m)$, let $s:=s_i$
denote the
least common multiple of $c_i$ and $c_{i+1}$ (resp., $c_i$ and $c_{i-1}$).
We say that $(c_i,m)$ is {\it significant} if $ms/c_i \in I_z$.
To $(c_i,m)$ we associate
the vector $\varphi=\varphi_{(c_i,m)}$
in $T_z \sn \simeq Hom_{k[X,Y]}(I_z,k[X,Y]/I_z)$
defined by
$$
\begin{array}{ccc}

\left \{
\begin{array}{l}
\varphi(c_j) = mc_j/c_i \mbox{ if } j \leq i \\
\varphi(c_j) = 0 \mbox{ if } j > i
\end{array}
\right .

&

\mbox{ resp. }

&

\left \{
\begin{array}{l}
\varphi(c_j) = 0 \mbox{ if } j < i \\
\varphi(c_j) = mc_j/c_i \mbox{ if } j \geq i
\end{array}
\right .

\end{array}
$$

{\fl According} to \cite[Theorem 3]{irreductible}, the set of elements
$\varphi_{(c,m)}$ for all significant cleft pairs $(c,m)$
is a basis of $T_z \sn$.

On the other hand, Nakajima gives a combinatorial description of the
weights occurring in $T_z \sn$ \cite[Proposition 5.8]{nakajima}.
Given a staircase $E$ and $e \in E$, let $a(e) := \max \{i \ |\ X^i
\cdot e \in E \}$ and let $b(e) := \max \{j \ |\ Y^j \cdot e \in E
\}$. The set of positive weights in $T_z \sn$, counted with
multiplicities, is the set of weights of the form
$w_+(e):=a(e)U-(b(e)+1)V$, and the set of negative weights is the
set of weights the form $w_-(e):=-(a(e)+1)U+b(e)V$.

\lpara

We now give a bijection $h_+$ resp. $h_-$ between the staircase $E$ and
the set of positive resp. negative
significant cleft couples, preserving the weights, meaning
that $wt(\varphi_{(h_\pm(e))}) = w_\pm(e)$.
The bijection $h_+$ resp. $h_-$
is defined as follows: given $e \in E$, we denote
$c_1 = Y^b \cdot e$ resp. $c_1 = X^a \cdot e$,
where $b$ resp. $a$ is the minimal integer such that
$Y^b \cdot e \not \in E$ resp. $X^a \cdot e \not \in E$.
We denote
$c_2 = X^a \cdot e$ resp. $c_2 = Y^b \cdot e$,
where $a$ resp. $b$ is the maximal integer such that
$X^a \cdot e \in E$ resp. $Y^b \cdot e \in E$.
We let $i$ resp. $j$ be the maximal integer
such that $X^{-i} \cdot c_1 \not \in E$ resp. $Y^{-j} \cdot c_1 \not \in E$
(thus $X^{-i} \cdot c_1$ resp. $Y^{-j} \cdot c_1$ is a cleft).
Finally we set
$h_+(e) = (X^{-i} \cdot c_1 , X^{-i} \cdot c_2)$ resp.
$h_-(e) = (Y^{-j} \cdot c_1 , Y^{-j} \cdot c_2)$:
by construction of $(c_1,c_2)$
this is a significant cleft couple.

\begin{prop}
\label{pro:2-tgt} The map $h_+$ resp. $h_-$ is a bijection between
$E$ and the set of positive resp. negative significant cleft
couples, and we have $wt(\varphi_{(h_+(e))}) = w_+(e)$ and
$wt(\varphi_{(h_-(e))}) = w_-(e)$.
\end{prop}
\begin{proo}
By symmetry, we give the proof only in the positive case. 
With the notations before the proposition, we have
$wt(\varphi_{(h_+(e))}) = wt(c_2c_1^{-1})$, which is readily $w_+(e)$. Thus
to prove the proposition it suffices to describe the inverse of $h_+$. Given
a positive significant cleft couple $(c,m)$, let $i$ be the maximal integer
such that $X^i \cdot m \in E$. The inverse of $h_+$ maps
$(c,m)$ to the greatest common divisor of $X^i \cdot c$ and $X^i \cdot m$.
\end{proo}








\subsection{Computation of the tangent space of $S^{[n,n+1]}$}

Let $(s,b) \in (\inc)^T$  ($s$ and $b$ stand for small and big).
We denote by $E_s,E_b$ their staircases.
Let $q:\inc \to S^{[n+1]}$ denote the natural projection,
and let $dq$ denote its
differential at $(s,b)$. There is a natural exact sequence
$$
0 \to \ker dq \to T_{(s,b)} \inc \stackrel{dq}{\to} T_b S^{[n+1]}\ .
$$

The following is immediate:

\begin{lemm}
\label{lem:incidence}
The tangent space $T_{(s,b)} S^{[p,q]}$ of $S^{[p,q]}$ at $(s,b)$
is the set of couples of homomorphisms
$(\varphi,\psi) \in Hom_{k[X,Y]}(I_s,k[X,Y]/I_s) \times
Hom_{k[X,Y]}(I_b,k[X,Y]/I_b)$ such that for all $f \in I_b$, we have
\[
\varphi(f)= \psi(f) \mbox{ mod } I_s.
\]
\end{lemm}
\begin{proo}
This tangent space is included in the tangent space
$T_{(s,b)} (S^{[p]} \times S^{[q]})$, which
is the direct sum $Hom_{k[X,Y]}(I_s,k[X,Y]/I_s) \oplus Hom_{k[X,Y]}(I_b,k[X,Y]/I_b)$.
Consider the ring $k[\epsilon]$, where $\epsilon^2 = 0$.
Identifying the tangent bundle of $S^{[p]}$ resp. $S^{[q]}$ with the
$k[\epsilon]$-points of $S^{[p]}$ resp. $S^{[q]}$, we see that
$(\varphi,\psi) \in T_{(s,b)} S^{[p,q]}$ if and only if the restriction of
$\varphi$ to $I_b$ is equal to the quotient of $\psi$ modulo $I_s$.
\end{proo}

The following propositions \ref{pro:image-generique} and
\ref{pro:image-speciale} describe the infinitesimal deformations of
$b$ that admit a lift to a deformation of $(s,b)$. There are two
cases, depending on the geometry of the staircases involved. 

If $m$ is a monomial, we denote by $x(m)$ its exponent for the
variable $X$. In other words $x(X^aY^b)=a$. Similarly $y(X^aY^b)=b$.
We denote by $c_1,\ldots,c_l$ the clefts of $s$, and by $k$ the index such that
$c_k \in E_b$.

\begin{prop}
\label{pro:image-generique}
Assume that $y(c_{k-1}) > y(c_k) + 1$, resp.
$x(c_{k+1}) > x(c_k) + 1$. Then the positive, resp. negative
part of
$\im dq$ is the subspace of
$T_b S^{[n+1]}$ generated by those $\varphi_{(c,m)}$ with
$(c,m) \not = (Yc_k,X^{-1}c_i) ,i>k$, resp.
$(c,m) \not = (Xc_k,Y^{-1}c_i),i<k$.
\end{prop}
\begin{proo}
By symmetry, it suffices to describe
the positive part of $\im dq$ when $y(c_{k-1}) > y(c_k) + 1$.

{\fl The fact} that many cleft pairs of $s$ are also cleft pairs of $b$ is a
potential source of confusion. Consequently, given a pair $(c,m)$ of both $s$
and $b$,
we will use $\varphi^{n}_{(c,m)}$ resp. $\varphi^{n+1}_{(c,m)}$ to denote the
corresponding tangent vector in $T_s S^{[n]}$ resp. $T_b S^{[n+1]}$.
Moreover we will use the convention that if $m \not \in E_s$, then
$\varphi^n_{(c,m)}=0$.
Let $(c,m)$ be a positive significant cleft pair of $b$.

{\fl First} note that if $c=c_i$ with $i \not = k$, then
$\varphi^{n+1}_{(c,m)} \in \im dq$. In fact in this case
it is clear that
$(\varphi^{n}_{(c,m)},\varphi^{n+1}_{(c,m)})$ is a tangent vector of $\inc$
(recall our convention that $\varphi^n_{(c,m)}$ is $0$ if $m=c_k$).

{\fl Next} consider the case where $c=Xc_k$. We also see that
$\varphi^{n+1}_{(c,m)} \in \im dq$ since the pair
$(\varphi^{n}_{(c_k,X^{-1}m)},\varphi^{n+1}_{(Xc_k,m)})$ is a tangent vector
of $\inc$ by Lemma \ref{lem:incidence}.

{\fl Now}, if $c=Yc_k$, the fact that $(c,m)$ is a significant cleft pair
implies that $Xm \in I_b$. If $m \not = X^{-1}c_i$ for any $i>k$, then
$XY^{-1}m \in I_b$ and thus $\varphi^n_{(c_k,Y^{-1}m)}(Xc_k)=0$, so that
$(\varphi^{n}_{(c_k,Y^{-1}m)},\varphi^{n+1}_{(Yc_k,m)})$ is a tangent vector
of $\inc$.

{\fl It remains} to show that the
$\varphi^{n+1}_{(Yc_k,X^{-1}c_i)}$-coefficient
of any vector in $\im dq$
vanishes, for all $i>k$. To this end, let
$(\varphi^n,\varphi^{n+1})$ be a tangent vector to the incidence variety.
Considering
$\varphi^n$ resp. $\varphi^{n+1}$
as an element of $Hom_{k[X,Y]}(I_s,k[X,Y]/I_s)$
resp. $Hom_{k[X,Y]}(I_b,k[X,Y]/I_b)$, we see
that the coefficient of $Y^{-1}c_i$ in $\varphi^n(Xc_k)$ is equal to the
coefficient of $X^{-1}c_i$ in $\varphi^n(Yc_k)$
(namely, those coefficients equal
the coefficient of $X^{-1}Y^{-1}c_i$ in $\varphi^n(c_k)$). It follows from
Lemma
\ref{lem:incidence} that
$\varphi^{n+1}$ has the same property. On the other hand, among all
the basis vectors $\left ( \varphi^{n+1}_{(c,m)} \right )$,
$\varphi^{n+1}_{(Yc_k,X^{-1}c_i)}$ is the only
vector for which these coefficients
are not equal. Thus the $\varphi^{n+1}_{(Yc_k,X^{-1}c_i)}$-coefficient in
$\varphi^{n+1}$ vanishes.
\end{proo}

\begin{prop}
\label{pro:image-speciale}
Assume that $y(c_{k-1}) = y(c_k) + 1$ resp.
$x(c_{k+1}) = x(c_k) + 1$. Then the positive resp. negative
part of
$\im dq$ is the subspace of
$T_b S^{[n+1]}$ generated by those $\varphi_{(c,m)}$ with
$(c,m) \not = (c_{k-1} , X^{x(c_{k-1})-x(c_k)-1} \cdot c_i),i>k$ resp.
$(c,m) \not = (c_{k+1} , Y^{y(c_{k+1})-y(c_k)-1} \cdot c_i),i<k$.
\end{prop}
\begin{proo}
The argument is similar to that used in the proof of the preceding proposition.
For any positive significant cleft pair $(c_i,m)$ with $i>k$ or $i<k-1$,
the vector
$( \varphi^n_{(c_i,m)} , \varphi^{n+1}_{(c_i,m)} )$ belongs to
$T_{(s,b)} \inc$; whence, $\varphi^{n+1}_{(c,m)}$ belongs to
the image of $dq$.

{\fl Similarly,} if $(c,m)$ is a positive cleft pair and $c=Xc_k$, then the
pair $(\varphi^n_{(c_k,X^{-1}m)}, \varphi^{n+1}_{(c,m)})$ belongs to
$T_{(s,b)} \inc$; we deduce
$\varphi^{n+1}_{(c,m)} \in \im dq$.

{\fl We} now consider positive significant cleft pairs of the form
$(c_{k-1},m)$. Assume that $i \geq k$
is such that
$y(c_{i-1})>y(m)\geq y(c_i)$.
Since
$(c_{k-1},m)$ is significant, we have
\[
x(m) \geq x(c_i) + x(c_{k-1}) - x(c_k) - 1.
\]
When $x(m) \geq x(c_i) + x(c_{k-1}) - x(c_k)$, the pair
$( \varphi^n_{(c_{k-1},m)} , \varphi^{n+1}_{(c_{k-1},m)} )$ belongs to
the tangent space
$T_{(s,b)} \inc$, so $\varphi^{n+1}_{(c_{k-1},m)} \in \im dq$. When
$x(m) = x(c_i) + x(c_{k-1}) - x(c_k) - 1$ and $y(m)>y(c_i)$, the pair
$( \varphi^n_{(c_k,X^{x(c_k)-x(c_{k-1})}Y^{-1}m)} ,
\varphi^{n+1}_{(c_{k-1},m)} )$
belongs to
$T_{(s,b)} \inc$, so $\varphi^{n+1}_{(c_{k-1},m)} \in \im dq$.

{\fl Finally} it remains to show that the
$\varphi^{n+1}_{(c_{k-1},X^{x(c_{k-1})-x(c_k)-1}c_i)}$-coordinate of
any vector in $\im dq$
vanishes. To this end, note that
$\varphi^{n+1}_{(c_{k-1},X^{x(c_{k-1})-x(c_k)-1}c_i)}$
is the only vector $\varphi$
in our basis of $T_b S^{[n+1]}$ for which
\[
[Y^{-1}c_i] \varphi(Xc_k) \neq [X^{-1}c_i] \varphi(Yc_k),
\]
while
\[
[Y^{-1}c_i]\psi(Xc_k) = [X^{-1}c_i]\psi(Yc_k)
\]
for any $\psi \in Hom(I_s,k[X,Y]/I_s)$, where $[m]\varphi(n)$ denotes
the coefficient of $m$ in $\varphi(n)$. 
\end{proo}

{\fl While} the description of $\im(dq)$ given in
Propositions ~\ref{pro:image-generique} and ~\ref{pro:image-speciale}
depends on the shape of the Young tableau corresponding to $s$,
the weights of the $T$-representation $T_b S^{[n+1]}/\im dq$
have a more uniform description. This space measures the obstructions
to lift a deformation of $b$ to a deformation of $(s,b)$. 

\beg{figure}

\end{figure}

\begin{prop}
\label{pro:coimage}
The weights of the $T$-representation $T_b S^{[n+1]}/\im dq$ have multiplicity
one and are given by
\[
(x(c_i)-x(c_k)-1) U + (y(c_i)-y(c_k)-1) V , i \not = k.
\]
\end{prop}
These weights are the weights of $\frac{c_i}{XYc_k}\ ,\ i \not = k$.
\begin{proo}
By symmetry, it suffices to consider only the positive part of the
quotient $T_b S^{[n+1]}/\im dq$.
Assume first that $y(c_{k-1}) > y(c_k) + 1$. By Proposition
\ref{pro:image-generique}, the positive weights
of $T_b S^{[n+1]}/\im dq$ are the
weights of $\frac{X^{-1}c_i}{Yc_k}$.
Thus the proposition is proved in this case.
Assume now that $y(c_{k-1}) = y(c_k) + 1$. By Proposition
\ref{pro:image-speciale}, the positive weights
of $T_b S^{[n+1]}/\im dq$ are the
weights of $\frac{X^{x(c_{k-1})-x(c_k)-1} c_i}{c_{k-1}}$.
Since $c_{k-1} = X^{x(c_{k-1})-x(c_{k})}Yc_k$, the proposition follows
in this case too.
\end{proo}

\begin{prop}
\label{pro:fibre}
The weights of $\ker dq$ have multiplicity one and are the following:
$$
\begin{array}{cc}
-(y(c_{i-1})-y(c_k)+1) V + (x(c_i)-x(c_k)-1) U & i>k \\
-(x(c_k)-x(c_{i+1})+1) U + (y(c_i)-y(c_k)-1) V & i<k
\end{array}
$$
\end{prop}
These weights are the weights of the arrows from $c_k$ to the corners of
the partition corresponding to $s$.
\begin{proo}
This kernel consists of
those morphisms
$\varphi \in Hom_{k[X,Y]}(I_s,k[X,Y]/I_s)$ for which $\varphi(I_b)=0$,
by Lemma \ref{lem:incidence}.
If follows that if $(c,m)$ is a cleft pair with $c \not = c_k$, we have
\[
[\varphi^n_{(c,m)}] \varphi = 0
\]
for every $\varphi \in \ker dq$. On the other hand, if $(c_k,m)$ is a
cleft pair and $\varphi$ has a non-vanishing
$\varphi^n_{(c_k,m)}$-coefficient, the fact that
\[
\varphi(Xc_k)=X\varphi(c_k)=0 \mbox{ and } \varphi(Yc_k)=Y\varphi(c_k)=0
\mbox{ mod } I_s
\]
implies that $m$ is a
corner of $s$. Conversely, for $m$ a corner of $s$, it is clear that
$\varphi^n_{(c_k,m)} \in \ker dq$. Thus $\ker dq$ is generated by
those elements $\varphi^n_{(c_k,m)}$ for which $m$ is a corner of $s$.
The proposition now follows immediately.
\end{proo}

The last two propositions and the exact sequence together describe the
tangent space $T_{(s,b)} \inc$ as a linear $T$-representation.
This will be useful later on to
compute equivariant Chern classes. Considering only the dimensions of
these spaces, we recover the following well-known result of Cheah
\cite[Theorem 3.2.2]{cheah} about the smoothness of $\inc$:

\begin{prop}
\label{pro:q1-lisse}
The incidence $\inc$ is a smooth irreducible subvariety of
$S^{[n]} \times S^{[n+1]}$.
\end{prop}
\comm
{
\begin{proo}
The incidence
 $\inc$ is connected of dimension $2n+2$ since any component contains
 a pair $(s,b)$ with $s$ and $b$ $T$-invariant and the smoothable
 component contains every such $(s,b)$ (this can be seen for example
using the family described in the proof of Proposition
\ref{pro:diviseur}).
On the other hand, by
Propositions \ref{pro:coimage} and \ref{pro:fibre}, the cokernel of $dq$ and
its kernel have the same dimension, namely, the number of clefts of $E_s$
minus one. Thus the tangent space of $\inc$
at $(s,b)$ has dimension $2n+2$. Since this is the dimension of
the smoothable component,
$(s,b)$
is a smooth point of $\inc$. Since the singular locus of $\inc$ is closed,
$T$-invariant, and has no $T$-fixed point, it is empty.
\end{proo}
}

\subsection{Application to the irreducibility of  $\inco$}

We consider the subvariety $S^{[p,q]}_0$ of $S^{[p,q]}$ parameterizing incident
schemes $(s \subset b)$ with respective length $p$ and $q$ both supported at the origin.
Recall Brian\c con's theorem which asserts the irreducibility of the variety $S_0^{[n]}\subset S^{[n]}$ parameterizing 
the subschemes of length $n$ and support the origin. The corresponding 
theorem for pairs of incident schemes is not true, as shown by the following example. 

\begin{prop}
\label{pro:s24}
The scheme $S^{[2,4]}_0$ is not irreducible.
\end{prop}
\begin{proo}
As $S^{[2,4]}_0$ is a strict subscheme of the irreducible 4-dimensional product
$S^{[4]}_0\times S^{[2]}_0$, 
any component has dimension 
at most 3 and any irreducible $F\subset S^{[2,4]}_0$  of dimension $3$ is an
irreducible component. One irreducible component 
of dimension 3 is birational to $S^{[4]}_0$ : its generic point parameterizes 
the couple $(s,b)$ with $b$ the generic curvilinear subscheme of length $4$
and $s$ the unique subscheme of 
$b$ with length $2$. An other 3-dimensional family is constructed as follows.
Let $(f,g)$ be two distinct linear forms and $b$ the subscheme of $S$ with
equation $(f^2,g^2)$. Let $s$ be any scheme of length $2$ supported at the
origin. Since $s\subset b$, the set of such $(s,b)$ describe a subvariety
$F\subset S^{[2,4]}_0$ of dimension $1+2=3$.
\end{proo}

Although the general incidence $S^{[p,q]}_0$ is wild and difficult to describe,
the case $(p,q)=(n,n+1)$ behaves nicely. 
The rest of this section is devoted to the proof of the following theorem: 
\begin{theo}
\label{theo:inco}
The incidence $\inco$ is irreducible of dimension $n$. At the
generic point $(s,b)$, the subschemes $s$ and $b$ are
curvilinear.
\end{theo}

We start the proof with a weaker version of the theorem,
in the next proposition. 

\begin{prop}
\label{prop-inco} We have $\dim \inco = n$ and there is only one
irreducible component $H_0$ in $\inco$ of dimension $n$. At a
generic point $(s,b) \in H_0$, the subschemes $s$ and $b$ are
curvilinear.
\end{prop}
\begin{proo}
Once again this result is a consequence of the detailed study of $\inco$
performed by Cheah \cite[Proposition 3.4.11]{cheah}. We give a short
proof.

To apply Bialynicki-Birula's decomposition theorem in $\inc$ which
is not compact, we first compactify $\inc$. So consider the
inclusion $\inc \subset \cinc$. Since $S^{[i,i+1]}$ is smooth for
all values of $i$, $\cinc$ is smooth. Consider the action of $k^*$
induced by the action on the affine plane defined by $t \cdot X =
t^\alpha X\ , \ t \cdot Y = tY$, where $\alpha$ is any integer
strictly greater than $n+1$.

Let $O$ denote the origin of $\A^2 \subset \p^2$. Let $Z$ be a
subscheme of $\p^2$. If $Z$ is not supported on $O$, then the limit
at $t=0$ of $t \cdot Z$ is also not supported on $O$. For $(s,b)$ a
$k^*$-fixed point in $\cinc$ let $C_{(s,b)}$ be the corresponding
Bialynicki-Birula cell: we have
$$
\displaystyle {
 \cinco = \coprod_{(s,b) \in { \left ( \inco \right ) }^T} C_{(s,b)}\ .
 }
$$
Thus to prove the proposition it is enough to show that all the
cells $C_{(s,b)}$ with $s$ and $b$ supported at the origin have
dimension at most $n$ and that exactly one has dimension $n$.

Let $(s,b) \in \inco$ be a given $k^*$-fixed point. Let us say that
a tangent vector $x \in T_{(s,b)} \inc$ is contractant if it is an
eigenvector for the $k^*$-action of positive weight. It is
well-known that the dimension of $C_{(s,b)}$ is the number of
independent contractant tangent vectors.

Let $x\in T_{(s,b)} \inc$ be a $T$-weight vector which is
contractant and let $w=aU+bV$ be its
weight. Its weight for $k^*$ is $a\alpha + b$ and this is a
positive integer. Since $\alpha>n+1$, we have $a \geq 0$ so $w$ is a positive
weight (recall the definition of positive weights in Subsection
\ref{subsection:tangent}). In particular, the vector space $W$ 
generated by such contractant tangent vectors $x$ satisfies $\dim W\leq
\dim T^+=n+1$, where $T^+\subset  T_{(s,b)} \inc$ is the vector space
generated by tangent
vectors with positive weight. Recall the description of $T_{(s,b)}
\inc$ given in terms of the projection $q:\inc \to S^{[n+1]}$ and
its differential $dq$. By Proposition \ref{pro:coimage}, all the
tangent vectors of $T_b S^{[n+1]}$ of positive weight $-V$ are in the image
of $dq$ but are not contractant. Since there is at least one such
vector, we get $\dim W\leq \dim T^+-1 = n$.
Moreover, if we have equality,
there is exactly one vector in  $T_b S^{[n+1]}$ of weight $-V$. 
This implies that the partition corresponding to $b$ is a rectangle:
$\lambda=(m,m,\ldots,m)$. But, if $m>1$, the eigenspace of weight $-V$
in $\im dq$ has dimension 1, as well as the eigenspace of weight
$-(m-1)V$ in $\ker dq$, by Proposition \ref{pro:fibre}. 
Since these
vectors are not contractant, we get $\dim W \leq n-1$ and a
contradiction. Thus the only possibility is $\lambda=(1,1,\ldots,1)$, and
the proposition is proved.
\end{proo}

The next proposition describes the Bialynicki-Birula cells of dimension 
$n-1$ introduced in the proof of the previous proposition.

\begin{prop} \label{prop:cellulesDeBBdeCodim1}
  A Bialynicki-Birula cell $C_{(s,b)}$ has dimension $n-1$ if and only if the
height of the partition $\lambda$ associated to $b$ is $2$. 
\end{prop}
\begin{proo}
We keep the notations of the previous proposition. 
  If the height $h$ of $\lambda $ is $1$, then $\dim C_{(s,b)}=n$ by the above. 
If $h\geq 3$, since there are at least $h$ independent
vectors in $\Im(dq)$ of weight a positive multiple of $-V$, 
there are at least three tangent vectors in $T_{(s,b)} \inc$
which are positive non contractant, hence
$\dim C_{(s,b)}\leq \dim T^+ - h\leq n-2$.
If $h=2$, then $\lambda=(2^\alpha,1^\beta)$ and the partition $\mu$ of
$s$ is $\mu_0=(2^{\alpha-1},1^{\beta+1})$  
or $\mu_1=(2^{\alpha},1^{\beta-1})$. Let us denote by $T_{cont}$ the
subspace of a vector space $T$  
generated by the contractant tangent vectors. Then
$\dim T_bS^{[n+1]}_{cont}=n-1$.

Our description of the kernel $\ker(dq)$ and the coimage $\coim(dq)$
of $dq$ (Propositions  \ref{pro:image-speciale} and \ref{pro:fibre})
show that when $\mu=\mu_1$ or $\beta=0$,
then $\dim \coim(dq)_{cont}=\dim \ker(dq)_{cont}=0$. When  $\mu=\mu_0$
and $\beta\neq 0$,
then $\dim \coim(dq)_{cont}=\dim \ker(dq)_{cont}=1$.  Summing up,
$\dim C_{(s,b)}= \dim (T_b)_{cont}-\dim \coim(dq)_{cont}+\dim
\ker(dq)_{cont}=n-1$  as required. 
\end{proo}

\begin{prop}\label{prop:dimDeLn}
  Let $L_{n}\subset \inc$ the set of $(s,b)$ such that $s$ is punctual. Then every component of 
$L_{n}$ has dimension at least $n+3$. 
\end{prop}
\begin{proo}

Following Gaffney and Lazarsfeld, if $f:X\rightarrow Y$ is a finite morphism between irreducible varieties 
we define the ramification locus $R_l\subset X$ containing the points $x$ for which $f^{-1}(f(x))$ is a scheme whose support on $x$ has length at least $l+1$. When $X$ is normal, $Y$ non-singular and $f$ surjective, 
then the components of $R_l$ have codimension 
at most $l$ \cite[p.58]{gaffneyLazarsfeld},\cite{lazarsfeldPHD}.

We apply this theorem with $X=U_n$  the universal family 
over $Y=\inc$ whose fiber over $(s,b)$ is the scheme $s$.
It suffices to prove that $U_n$ is normal. We shall prove that $U_n$ is  
Cohen-Macaulay and smooth in codimension one, which implies normality according to Serre's criteria. 

$U_{n}$ is Cohen-Macaulay as it is flat over the smooth base $\inc$. 

For any $\lambda=(\lambda_1,\dots,\lambda_k)$ ordered $k$-tuple with 
$\sum \lambda_i=n+1$, we denote by $\Delta_{\lambda}\subset S^{[n+1]}$ the 
stratum of subschemes $z$ of type $\lambda$, ie.
$z=z_1\amalg \dots \amalg z_k$ with
$length(z_k)=\lambda_k$ and $z_k$ punctual. Since any punctual
$z_i(p)$ supported by $p$ is the translation of a subscheme $z_i(0)$
supported by the origin, $\dim \Delta_{\lambda}=\dim 
(S^k\times S^{[\lambda_1]}_{0}\times  \dots \times S^{[\lambda_k]}_{0})$
 
For $i\leq k$, let
$\mu_i(\lambda)=(\lambda_1,\dots,\lambda_{i-1},\lambda_i-1,\lambda_{i+1},\dots)$. 
For $\mu \subset \lambda$ with $\Sigma \mu_i = n$ let
$D_{\mu,\lambda}\subset \inc$ be the image
of the generically well defined quasi-finite map 
\begin{eqnarray*}
S^k \times S_0^{[\mu_1,\lambda_1]}\times 
\dots \times S^{[\mu_k,\lambda_k]}_0 &\rightarrow &\inc\\
((p_1,\dots,p_k),(t_1,w_1),\dots,(t_k,w_k))& \mapsto &
(\amalg t_i(p_i),\amalg w_i(p_i)).
\end{eqnarray*}

Let $D_\lambda\subset \inc$ denote the inverse image of
$\Delta_\lambda$ by the natural projection $\inc \rightarrow S^{[n+1]}$. 
Then $D_\lambda=\cup_{i\leq k}D_{\mu_i(\lambda),\lambda}$. 
For $\lambda\neq (2,1,\dots1)$ and $\lambda\neq
(1,\dots,1)$, the codimension of  $D_\lambda$ in $\inc$ is at least
$2$  according to Proposition \ref{prop-inco}.
In particular, 
no smoothness condition is required for the universal family $U_n$ over
$D_\lambda$. When $\lambda=(1,\ldots,1)$, the smoothness of $U_n$ is obvious.

We consider now the case $\lambda= (2,1,\dots1)$.
Let $(s,b)\in D_\lambda$.

If $(s,b)\in D_{\mu,\lambda}$ with
$\mu=(1,\dots,1)$, then 
locally around $(s,b)$, $\inc$ is isomorphic to 
$S^{[1,2]}\times S^{n-1}$. The universal family $U_n$ over
$D_{\mu,\lambda}$ is locally a
disjoint union of sheets. The sheets coming from the universal
families over $S$ are obviously smooth. The last sheet $Z$ coming from the 
factor $S^{[1,2]}$ is such that the projection $Z\fd S^{[1,2]}$ is an
isomorphism (the fiber is zero dimensional with length 1), 
so this last sheet is smooth too.

If $(s,b)\in D_{\mu,\lambda}$ with
$\mu=(2,1,\dots,1,0,1,\dots,1)$, then 
locally around $(s,b)$, $\inc$ is isomorphic to $S^{[2]}\times S^{n-1}$.
The universal family $U_{n}$ is smooth since it is the disjoint union of the
pullback of the smooth universal families over $S^{[2]}$ and $S$.

\end{proo}

\begin{coro}
\label{coro:dimPossibleDesComposantes}
  Every irreducible component of $\inco$ has dimension $n$ or $n-1$. 
\end{coro}
\begin{proo}
  Moving the subschemes of $\inco$ with translations, the product
  $\inco \times S$ parameterizes the set $L$ of pairs $(s,b)\in \inc$
  with $s,b$ punctual with any support $p\in S$. We need to prove that the
  components of $L$ have dimension $n+1$ or $n+2$. Consider the
  residual morphism $Res:\inc \rightarrow S$ that sends a pair $(s,b)$ to
  the point $q$ defined by the ideal $(I_s:I_b)$. Let $\Delta_x:L_n \fd k$,
$(s,b)\rightarrow x(s)-x(Res(s,b))$,
where $x(s)$ denotes the $x$ coordinate of the punctual subscheme $s$.
Define similarly $\Delta_y$.
The components of $L_n$ have dimension at least $n+3$ by
Proposition \ref{prop:dimDeLn}. 
From the equality $L=L_n \cap \Delta_x^{-1}(0)
\cap \Delta_y^{-1}(0)$, we conclude that any component of 
$L$ has dimension at least $n+1$.
The components of $L$ have dimension at most $n+2$ by
Proposition \ref{prop-inco}. 
\end{proo}

We now conclude the proof of Theorem \ref{theo:inco}. For the proof,
we need to produce some universal families over Bialynicki-Birula
cells. They are constructed from
the description of the tangent space by a procedure similar to the one used in 
\cite{irreductible}.
\begin{proo}
By Corollary \ref{coro:dimPossibleDesComposantes},
the components of $\inco$ have dimension $n$ or $n-1$.
The Bialynicki-Birula decomposition of $\inco$ is a partition into irreducible
sets. It follows that 
the irreducible components of $\inco$ are the maximal sets for the inclusion
among the closure of 
the Bialynicki-Birula cells. Since we already proved that there is a unique
maximal component 
of dimension $n$, it remains to prove that the closure of the cells of
dimension $n-1$ described by Proposition 
\ref{prop:cellulesDeBBdeCodim1} are not irreducible components of $\inco$. 

Let $C_{s_0,b_0} \subset \inco$ be a Bialynicki-Birula
cell of dimension $n-1$, $(s_0,b_0)\in \inco$ the corresponding fixed point,
$\lambda=(2^\alpha,1^\beta)$ and $\mu$ be the partitions of $b_0$ and $s_0$.
Let $\mu_0=(2^{\alpha-1},1^{\beta+1})$ 
and $\mu_1=(2^{\alpha},1^{\beta-1})$. We have $\mu=\mu_0$ or $\mu=\mu_1$. 

First, we remark that the irreducible components of $\inco$ are invariant under
$GL_2$, the group of linear automorphisms of the plane.
In particular, if we prove that
the generic
point of $C_{s_0,b_0}$ is not invariant under $GL_2$, it follows
that the closure $\overline{C_{s_0,b_0}}$ is not a component of $\inco$. 
Moreover we will use the following notation: we denote by
$k[X,Y]_d$ the space of homogeneous polynomials of degree $d$ and by
$\pi_d:k[X,Y] \to k[X,Y]_d$ the natural projection. Given an ideal $I$,
$I_d$ will denote the subspace $\pi_d(I) \subset k[X,Y]_d$. Let us finally
say that an admissible cleft couple for $I_b$ is liftable if the corresponding
infinitesimal deformation of $I_b$ can be lifted to a infinitesimal
deformation of the pair $(I_s,I_b)$.

\lpara
$\bullet$ If $\alpha>1$ and $\beta=0$, then consider
the ideal $I_b=(X^{\alpha},Y^2+\sum_{0<i<\alpha,j\leq 1}c_{ij}X^i Y^j)$ and 
$I_s=I_b+(YX^{\alpha-1})$. The variables $c_{ij}$ are the $(n-1)$
coordinates on the Bialynicki-Birula cell $C$ and $I_s$,$I_b$ are the
corresponding universal ideals. If $(I_s,I_b)$ is the generic element in this
cell, we have ${(I_b)}_1 = k \cdot X$, thus this generic point is not
$GL_2$-invariant. 

\lpara
$\bullet$
If $\alpha>1$, $\beta=1$ and $\mu=\mu_1$, we again have 
${(I_b)}_1 = k \cdot X$ for the generic pair $(I_s,I_b)$, and this is not
$GL_2$-invariant.

\lpara
$\bullet$
If $\alpha=1$ and $\beta=0$, then $I_b=(X,Y^2)$ and $I_s=(X,Y)$, so this
point is not $GL_2$-invariant. Similarly, if
$\alpha=1,\beta=1$ and $\mu=\mu_1$, then
$I_s=(X,Y^2)$, and this is not $GL_2$-invariant.

\lpara
$\bullet$
If $\alpha > 1$ and $\beta >2$ then the cleft couple
$((2,0),(0,2))$ is not admissible for
$I_b$. It follows that the generic element $(b,s)$, ${(I_b)}_2$ has dimension one and is generated
by a polynomial divisible by $Y$. 
In particular, the generic element of the cell is not $GL_2$-invariant.

\lpara
$\bullet$
If $\alpha>1$, $\beta=2$ and $\mu=\mu_0$, then the cleft couple
$((2,0)(0,2))$ is admissible, but not
liftable, by
Proposition \ref{pro:coimage}. Thus for the same reason the cell is not
$GL_2$-invariant.

\comm{
$\bullet$ If $\beta>1$ and $\alpha>1$, then one may similarly write down the
ideals $I_b,I_s$ of the universal families over the Bialynicki-Birula cell.
If $\mu=\mu_0$,
$I_b=(X^{\alpha+\beta},YX^{\alpha}+\sum_{\alpha<j<\alpha + \beta}
a_jX^j,Y^2+\sum_{0<j<\alpha}b_jYX^j+\sum_{\beta<j<\alpha+\beta}c_jX^j)$
and $I_s=I_b+(YX^{\alpha-1}+dX^{\alpha+\beta-1})$. Note that there is
no coefficient $c_j$ for $j=\beta$ since 
 $((0,2),(\beta,0))$ is a cleft couple corresponding to a 
vector which is not in $\im dq$. If $\mu=\mu_1$,
$I_b=(X^{\alpha+\beta},YX^{\alpha}+\sum_{\alpha<j<\alpha+\beta}a_jX^j,Y^2+
\sum_{0<j<\alpha}b_jYX^j+\sum_{\beta\leq j<\alpha+\beta}c_jX^j)$
and $I_s=I_b+(X^{\alpha+\beta-1})$. In particular the degree $2$ part
of any element of $I_b$ contains no term $X^2$ but may contain a term
$Y^2$ when $\mu=\mu_0$ or when ($\mu=\mu_1$ and $\beta>2$). In
geometric terms, the singularity associated with the generic subscheme 
characterizes the tangent of the branch $Y=0$ and the generic point is
not invariant under
every linear transformation.
}

\lpara
It remains to consider the cases ($\beta=1$ and $\mu=\mu_0$),
$(\alpha=1,\beta>1,\mu=\mu_0)$, $(\alpha=1,\beta>1,\mu=\mu_1)$,
and $(\alpha>1,\beta=2$ and $\mu=\mu_1$).
For these cases, we will see that
the closure of the Bialynicki-Birula cells are invariant under 
$GL_2$, and thus we cannot apply the same 
arguments as above. Instead, we will prove that the closure 
$\overline {C_{s_0,b_0}}$ of the cell 
under consideration is not an irreducible component of $\inco$
as it is included in the unique (``curvilinear'') component of dimension $n$. 
To this end, we apply a change of coordinates to obtain simple equations
for the generic point $(s,b)$ of the cell $C_{s_0,b_o}$
and we express $(s,b)$ as the limit of
$(s(t),b(t))\in \inco$ with $s(t)$ and $b(t)$ curvilinear.

\lpara 
$\bullet$ Consider the case $\beta=1$ and $\mu=\mu_0$. The universal
families over $C_{s_0,b_0}$ are described by coordinates $c_{ij},d$ and
universal ideals
$I_b=(X^{\alpha+1},YX^\alpha,Y^2+\sum_{i+j\geq 2,(i,j)\in \lambda} c_{ij}X^iY^j)$,
$I_s=I_b+(YX^{\alpha-1}+d X^\alpha)$. Note that $I_b$ contains all the 
monomials of degree $\alpha+1$. The element in
$I_b$ with initial term $Y^2$ vanishes on a
curve locally reducible around $0$ as a union of two distinct smooth
curves when the coefficients are generic. Up to a change of
coordinates, one may suppose that the two branches have equations
$X=0$ and $Y=0$. 
Then $I_b$ contains $XY$ and all the monomials of degree $\alpha +1$.
Thus  $I_b=(X^{\alpha+1},Y^{\alpha+1},XY)$ because of the inclusion, 
and both ideals have the same
colength. The ideal $I_s$ has codimension one in $I_b$
thus the general element has the form
$I_s=I_b+(X^{\alpha}+dY^{\alpha})$. Up to a linear change of coordinates
of the form
$Y\mapsto c.Y$, $X\mapsto X$, one may suppose that
$I_s=I_b+(X^{\alpha}+Y^{\alpha})$.
For a generic pair $(s,b)$, $b$ is the union of two curvilinear schemes
of length $\alpha+1$ supported at the origin, and $s$ a colength one
subscheme of $b$: this cell is invariant under automorphisms.

For $t \not = 0$ let $I_b(t)=(XY+t(X+Y),X^{2\alpha+1},Y^{2\alpha+1})$ and
$I_s(t)=(XY+t(X+Y),X^{2\alpha},Y^{2\alpha})$.
Let $(I_s(0),I_b(0))$ be the limit at $t=0$ of $(I_s(t),I_b(t))$.
Obviously $XY=\lim_{t \rightarrow 0} XY+t(X+Y)\in I_b(0)$. 

Since $b(t)$ (resp. $s(t)$) has
length $2\alpha+1$ (resp. $2\alpha$) with support the
origin, every monomial of degree $2\alpha+1$ (resp. $2\alpha$)
is in $I_b(t)$ (resp $I_s(t)$).
Since $X+Y= \frac{-XY}{t}$ modulo $I_b(t)$, we obtain
$X(X+Y)^\alpha\in I_b(t)$ and
$Y(X+Y)^\alpha\in I_b(t)$. Summing up,
$I_b(0)\supset
(XY,X(X+Y)^\alpha,Y(X+Y)^\alpha)=(XY,X^{\alpha+1},Y^{\alpha+1})$. This
inclusion is an
equality since the two ideals have colength $2\alpha+1$.
The same reasoning with the curvilinear $s(t)$
instead of $b(t)$ 
shows that $I_s(0)\supset (XY,X^{\alpha+1},Y^{\alpha+1})$.
Modulo $I_s(t)$, $(X+Y)^\alpha=(\frac{-XY}{t})^\alpha=0$.
Thus $I_s(0)\supset  (XY,X^{\alpha+1},Y^{\alpha+1},(X+Y)^\alpha)$
and the equality follows by length
considerations. We have proved $I_b=I_b(0)$ and $I_s=I_s(0)$, as expected.  

\lpara
$\bullet $ If $\alpha>1,\beta=2,\mu=\mu_1$, we can perform as above a change of
coordinates in order to reduce to the case
\comm{
$I_b=(X^{2+\alpha},YX^\alpha+aX^{1+\alpha}, Y^2+
\sum_{1\leq i<\alpha}b_iYX^i +\sum_{1<i<\alpha+2}c_iX^i$,
$I_s=I_b+(X^{1+\alpha})$. 
Remark that the generic element of $I_b$ with
initial term $Y^2$ is a union of two smooth curves in the neighborhood
of zero. If we choose a system of coordinates such that this curve 
the two local branches are $Y=X$ and $Y=-X$, then the ideals are 
$I_b=(X^{2+\alpha},YX^\alpha+aX^{1+\alpha}, Y^2-X^2)$,
$I_s=I_b+(X^{1+\alpha})=(X^{1+\alpha}, YX^\alpha,Y^2-X^2)$. 
Thus $I_s$ defines the union of the two curvilinear subschemes of $(X^{1+\alpha},
Y-X)$ and $(X^{1+\alpha},
Y+X)$. Since $I_b\supset (X^{2+\alpha},Y^2-X^2)$, the associated
subscheme $Z(I_b)$ is included in the union 
of the two curvilinear subschemes  $(X^{2+\alpha},
Y-X)$ and $(X^{2+\alpha},Y+X)$. We make a further change of coordinate 
so that the branches $X+Y=0$ and $X-Y=0$ become $X=0$ and $Y=0$.
}
$I_s=(XY,X^{\alpha+1},Y^{\alpha+1})$ and
$I_b=((XY,X^{\alpha+2},Y^{\alpha+2},X^{\alpha+1}+Y^{\alpha+1})$.
For a generic pair $(b,s)$, $b$ is a colength one subscheme in
the union $c_1 \cup c_2$ of two curvilinear subschemes of length
$\alpha +2$, and $s$ is the union $c_1' \cup c_2'$, where $c_i' \subset c_i$
is the unique colength one subscheme: this cell is invariant under
automorphisms.
The same computation as above
now shows that $I_s=\lim I_s(t)$, $I_b=\lim I_b(t)$ with 
$I_s(t)=(XY+t(X+Y),X^{2\alpha+1},Y^{2\alpha+1})$ and
$I_b(t)= (XY+t(X+Y),X^{2\alpha+2},Y^{2\alpha+2})$.

\lpara
$\bullet $ If $\alpha=1,\beta>1,\mu=\mu_0$, then
$I_b=(X^{1+\beta},XY+\sum_{2\leq j \leq \beta} a_jX^j, Y^2+\sum_{2\leq
  j \leq \beta} a_jYX^{j-1})$ and
$I_s=I_b + (Y+\sum_{2\leq j \leq \beta} a_jX^{j-1}+dX^\beta)$. 
Up to the
coordinate change 
$X\mapsto  X$, $Y\mapsto Y+\sum_{2\leq j \leq \beta} a_jX^{j-1}+dX^\beta$,
one may suppose that
$I_b=(X^{1+\beta},XY,Y^2) $ and $I_s=(X^{1+\beta},Y)$.
It follows that for a generic pair $(b,s)$, $s$ is a curvilinear scheme
and $b$ the union of $s$ and the $2$-fat point: this cell is invariant under
automorphisms. 
Consider the
curvilinear ideals $c=(X^{2+\beta},Y)$,  and the automorphism
$\phi_t:X\mapsto X, Y\mapsto tY+X^{\beta+1}$. The ideals
$I_b(t)=\phi_t(c)$ and $I_s(t)=\phi_t(I_s)=I_s$ are such that
$\lim_{t\rightarrow 0}I_s(t)=I_s$ and $\lim_{t\rightarrow 0}I_b(t)=I_b$.

\lpara
$\bullet$
If $\alpha=1,\beta>1,\mu=\mu_1$, then
$I_b=(X^{1+\beta},XY+\sum_{2\leq j \leq \beta} a_jX^j, Y^2
+\sum_{2\leq j \leq \beta} a_jYX^{j-1}+dX^\beta)$ and $I_s=I_b+(X^\beta)$. 
Up to the 
two coordinate changes
$X\mapsto  X$, $Y\mapsto Y+\sum_{2\leq j \leq \beta} a_jX^{j-1}$,
and then $X\mapsto X, Y\mapsto \lambda Y$, one may suppose that
$I_b=(X^{1+\beta},XY, Y^2+X^\beta)$ and $I_s=(X^{\beta},XY, Y^2)$.

For a generic pair $(b,s)$, $b$ is a colength one subscheme in
the union $c_1 \cup c_2$ of two curvilinear subschemes of length
$3$ and $n-1$, and $s$ is the union $c_1' \cup c_2'$, where $c_i' \subset c_i$
is the unique colength one subscheme: this cell is invariant under
automorphisms.

For $t \not = 0$ let $I_b(t)=(XY-t^2Y+t^\beta X,X^{\beta+2},Y^{\beta+2})$ and
$I_s(t)=(XY-t^2Y+t^\beta X,X^{\beta+1},Y^{\beta+1})$.
Let $(I_s(0),I_b(0))$ be the limit at $t=0$ of $(I_s(t),I_b(t))$.
Obviously $XY=\lim_{t \rightarrow 0} XY -t^2Y+t^\beta X \in I_b(0)\cap I_s(0)$. 

Since $b(t)$ (resp. $s(t)$) has
length $\beta+2$ (resp. $\beta+1$) with support the
origin, all the monomials in $k[X,Y]$ of degree $\beta+2$ (resp. $\beta+1$)
are in $I_b(t)$ (resp $I_s(t)$). A straigthforward induction shows
that
\begin{displaymath}
  \forall k\geq 1,\
  Y=t^{\beta-2}X+t^{\beta-4}X^2+\dots+t^{\beta-2k}X^k+t^{-2k}X^kY \ mod\
  I_b(t) \cap I_s(t). 
\end{displaymath}
In particular,
$e(t):= Y-t^{\beta-2}X-\cdots -t^{-\beta-2}X^{\beta+1} \in I_b(t)$,
$f(t):=Y-t^{\beta-2}X-\cdots -t^{-\beta}X^\beta \in I_s(t)$, 
and 
\begin{eqnarray*}
X^{\beta+1}=\lim_{t\rightarrow 0}t^{\beta+2}e(t)\in I_b(0),\\
X^\beta=\lim_{t\rightarrow 0}t^{\beta}f(t)\in I_s(0).
\end{eqnarray*}
Since $Y^2=(Y-e(t))^2\ \ mod\ I_b(t)$ and $Y^2=(Y-f(t))^2\ \ mod\ I_s(t)$,
we get 
\begin{eqnarray*}
  g(t):=Y^2-t^{2\beta-4}X^2-\cdots-(\beta-1)X^\beta-\beta t^{-2} X^{\beta+1}\
  \in I_b(t),\\
 h(t):=Y^2-t^{2\beta-4}X^2-\cdots-(\beta-1)X^\beta\
  \in I_s(t).
\end{eqnarray*}
It follows that 
\begin{eqnarray*}
Y^2+X^\beta=\lim_{t\rightarrow 0}g(t)-\beta t^\beta e(t) \in I_b(0)\\
Y^2=\lim_{t\rightarrow 0} h(t)-(\beta-1) t^\beta f(t)  \in I_s(0)
\end{eqnarray*}
Summing up, these limits prove that $I_s(0)\supset I_s$ and $I_b(0)\supset I_b$. 
The equalities follows from the inclusions by length considerations. 
\end{proo}

\section{Bases of the equivariant Chow ring}

\label{section:bases}

We now present three natural bases
$\fix(\lambda),\nak(\lambda),\es(\lambda)$ of the $K$-vector space
$A^*_K(\sn)$.  Our three bases
of $\ck{\sn}$ are naturally parameterized by the set
$\pn$  of partitions $\lambda$ of
weight $n$.

\lpara

Let $n \geq 0$ and let $i > 0$ be integers. We define some
correspondences following Nakajima \cite{nakajima}:

\begin{defi}
\label{def:qi}
  Let $Q_i^n \subset \sn \times S^{[n+i]}$ be the
  closure of the set of pairs $(z_n,z_{n+i})$ where $z_n \in \sn$ is arbitrary
  and $z_{n+i} \in S^{[n+i]}$ is the disjoint union of $z_n$ and a punctual
  scheme of length $i$.
\end{defi}

The $T$-invariant correspondence $\amalg_n\ Q_i^n$ induces an
operator (called ``creation operator'') $\displaystyle{q_i:\oplus_n
\ A^*_T(\sn) \to \oplus_n \  A^{*+i-1}_T(S^{[n+i]})}$ on Chow
groups. Assume now that $i<0$. The ``destruction operator'' $q_i$ is
defined either as the dual of $q_{-i}$ or with the correspondence
$Q_i^n \subset \sn \times S^{[n+i]}$ which is dual to the
correspondence $Q_{-i}^{n+i}$, in the sense of Definition
\ref{def:corr-duale}. By Proposition \ref{pro:duale}, both
definitions lead to the same operator. For any $i$, $q_i$ has
conformal degree $i$ and cohomological degree $i-1$. We make the
convention that $q_0=0$.

Given a partition $\lambda$ of length $l$, we will denote by
$\nak(\lambda)$ the equivariant class obtained applying
$q_{\lambda_l} \circ \cdots \circ q_{\lambda_1}$ to the vacuum $\vacuum$,
where the vacuum is the fundamental class on $S^{[0]}$.
Since by \cite{nakajima}, the classes $\nak(\lambda)$ for $\lambda
\in \pn$ restrict to a basis of the non equivariant Chow group of
$\sn$, $\{\nak(\lambda) \ , \ \lambda \in \pn \}$ is a basis of
$A^*_K(\sn)$ over $K$.

\lpara

Recall the classes introduced by Ellingsrud and Str\o mme in
\cite{ES}. These classes are introduced for $\p^2$ but we can
consider the same classes for $\A^2$. We choose an injection
$k^* \to T,t\mapsto (t^{-1},t^{-d})$
where $d$ is large. The action of $T$ on $S^{[n]}$ induces an action of $k^*$. 
With the assumption that $d$ is large enough, the 
$k^*$-fixed points are the $T$-fixed points; in particular there is
a finite number of them and they are parameterized by partitions. 
More precisely, if $\lambda=(\lambda_1,\dots,\lambda_l)$ is a
partition, we denote by $x_{\lambda}$ the subscheme with ideal 
$I_{x_{\lambda}}$ generated by the $l+1$ polynomials
$X^{k-1}Y^{\lambda_{k}}$,
where $k$ varies from 1 to $l+1$, with the convention that
$\lambda_{l+1}=0$. 

To each partition $\lambda$ of weight $n$ corresponds a
Bialynicki-Birula cell containing the points $p\in S^{[n]}$ 
such that $\lim_{t \rightarrow 0}t.p=x_{\lambda}$.
We denote $ES_\lambda \subset S^{[n]}$
the closure of this cell. Let $l$ be the length of the partition
$\lambda$. 
Geometrically, the Bialynicki-Birula cell associated to $\lambda$ 
parameterizes the subschemes $Z\subset S$ for which
there exist $x_1,\ldots,x_l \in k$ such that each
intersection $Z \cap \{X=x_i\}$ has length $\lambda_i$. The
equivariant class of $ES_\lambda$
in the Chow ring will be denoted $\es_\lambda$. Since by
definition $S^{[n]}$ has the cellular decomposition $S^{[n]} =
\amalg_\lambda ES_\lambda$, where $\lambda \in \pn$, the classes
$\es_\lambda$ for $\lambda \in \pn$ form a basis of $\ch{\sn}$.

\lpara

Finally, the classes $\fix(\lambda) \in \ck{\sn}$ are defined using
the localization theorem \cite[Theorem 1]{EG}.
The set
${(\sn)}^T$ contains the points $x_{\lambda}$ parameterized by $\lambda
\in \pn$. Let $1_\lambda \in
A^*_T({(\sn)}^T)$ be the class corresponding to
$x_\lambda$. Let $i:{(\sn)}^T \to \sn$ denote the inclusion.
By Lemma \ref{lem:i*},
$i^*_K:A^*_K(\sn)^T \rightarrow A^*_K(\sn)$ is an isomorphism.

\begin{defi}
\label{def:fix}
Let $\fix(\lambda)$ be the unique element in
$\ck{\sn}$ such that $i^*_K(\fix(\lambda)) = 1_\lambda$.
\end{defi}

\noindent
Let us denote by $\Tan(\lambda) \in \Z[U,V]$ the product of the weights of the
tangent space $T_{x_\lambda}\sn$.
According to the self-intersection formula, we have
\begin{equation}
\label{equ:fix} \fix(\lambda) = i_*(1_\lambda) / Tan(\lambda) \ .
\end{equation}

\noindent
Recall Definition \ref{def:produit-scalaire}.
We deduce from (\ref{equ:fix}) the following lemma:

\begin{lemm}
\label{lem:norme-fix} We have
$\scal{\fix(\lambda),\fix(\lambda)}_\sn = 1/\Tan(\lambda)$.
\end{lemm}
\begin{proo}
By (\ref{equ:fix}), we have
$$
\scal{\fix(\lambda),\fix(\lambda)}_\sn = \frac
{\scal{i_*(1_\lambda),i_*(1_\lambda)}_\sn} {\Tan(\lambda)^2}.
$$
By Definition \ref{def:produit-scalaire}, this is
$\pi^K_*(i_*(1_\lambda) \cup i_*(1_\lambda))/\Tan(\lambda)^2$,
if $\pi:\sn \to \Spec k$
denotes the projection to a point. Since
$i_*(1_\lambda) \cup i_*(1_\lambda) = \Tan(\lambda)\cdot i_*(1_\lambda)$,
the lemma follows.
\end{proo}

\section{Classical Operators}

\label{sec:classical}

Let us denote by $A$ the direct sum $\bigoplus_n \ch{\sn}$
and $A_K:=\bigoplus_n \ck{\sn}$. In this
section, we consider the classical operators acting on $A_K$, namely
the creation/destruction operators $q_i$ and the boundary operator
$\partial$, and an auxiliary operator $\rho$. We compute them in the
basis $\fix(\lambda)$. We also compute the commutators of these
operators.

The operators $\partial, \rho,q_i$ for $i>0$ are naturally defined on
$A$ and they are naturally extended to $A_K$. We use freely the same 
notation for the operators on $A$ and on $A_K$. On the contrary, 
the
operators $q_i$ for $i<0$ are defined on $A_K$ but not on $A$. 
This is because their definition involves non proper morphisms. 

In Theorem \ref{theo:qi} we give an explicit algorithm to compute
all operators $q_i$ in the basis $\fix(\lambda)$. With the help of
this result, we checked on a computer our formulas for commutators,
such as Theorem \ref{theo:qi-qj}. However the computations are very
tricky and we are not able to give a purely algebraic proof of these
formulas: instead we use geometric arguments.

\subsection{The operators $q_1$ and $q_{-1}$}

Given a partition $\lambda$, we denote by $\lambda[1]$ the set of
partitions $\mu$
with $\lambda \subset \mu$ and
$|\mu| = |\lambda| + 1$.
Given two partitions $\lambda,\mu$ with $\mu \in \lambda[1]$, we denote by
$\Cok(\lambda,\mu) \in \Z[U,V]$ the product of the weights of Proposition
\ref{pro:coimage} and by $\Ker(\lambda,\mu) \in \Z[U,V]$ the product
of the weights of Proposition \ref{pro:fibre}.
With these notations we have the following proposition:

\begin{prop}
\label{pro:q1}
We have the following formula:
$$
q_1(\fix(\lambda)) = \sum_{\mu \in \lambda[1]} \frac
{\Cok(\lambda,\mu)}{\Ker(\lambda,\mu)} \fix(\mu) \ .
$$
\end{prop}
\begin{proo}
Consider the following diagram:
$$
\xymatrix
{
& \sn \times S^{[n+1]} \ar[ld]_{\pi_{n}} \ar[rd]^{\pi_{n+1}} \\
S^{[n]} & & S^{[n+1]}
}
$$

\noindent
Let $[S^{[n,n+1]}] \in \ch{\sn \times S^{[n+1]}}$ denote the
equivariant class of the incidence. By definition, we have
$q_1(\fix(\lambda)) = \pi_{n+1,*}^K \left (
\pi_{n,K}^*(\fix(\lambda)) \cup [S^{[n,n+1]}] \right )$. Now, by
Definition \ref{def:fix} and the following commutative diagram
$$
\xymatrix
{
S^{[n]} \times S^{[n+1]} \ar[d] & {(S^{[n]})}^T \times {(S^{[n+1]})}^T
\ar[l] \ar[d] \\
S^{[n]} & {(S^{[n]})}^T\ ,  \ar[l]
}
$$
we have $\pi_{n,K}^* (\fix(\lambda)) = \sum_\mu \fix(\lambda)
\otimes \fix(\mu)$. By Proposition \ref{pro:q1-lisse}, $S^{[n,n+1]}$
is smooth, so the restriction of its class to a fixed point is the
product of the weights of the normal space at this point. By
Propositions \ref{pro:coimage} and \ref{pro:fibre}, we thus have
$$
[S^{[n,n+1]}] = \sum_{\lambda,\mu:\mu \in \lambda[1]}
\frac{\Tan(\lambda)\Cok(\lambda,\mu)}{\Ker(\lambda,\mu)}\ \
\fix(\lambda) \otimes \fix(\mu).
$$
Therefore, \vspace{-5mm}
$$
\begin{array}{rcl}
q_1(\fix(\lambda)) & = & \pi_{n+1,*}^K \left ( \displaystyle {
\sum_{\mu \in \lambda[1]} \frac
{\Tan(\lambda)\Cok(\lambda,\mu)}{\Ker(\lambda,\mu)}\ \
\fix(\lambda) \otimes \fix(\mu) }
\right ) \\
& = & \displaystyle { \sum_{\mu \in \lambda[1]} \frac
{\Cok(\lambda,\mu)}{\Ker(\lambda,\mu)} \fix(\mu)\ , }
\end{array}
$$
using Theorem \ref{theo:bott}.
This is what we wanted to prove.
\end{proo}

For example, we have $q_1(\fix([2])) = \frac{-2U+V}{-U+V}
\fix([2,1]) + 3 \fix([3])$. This is illustrated as follows, where
the weights of the blue resp. red arrows are the numerators resp.
denominators of
the coefficients:

$$
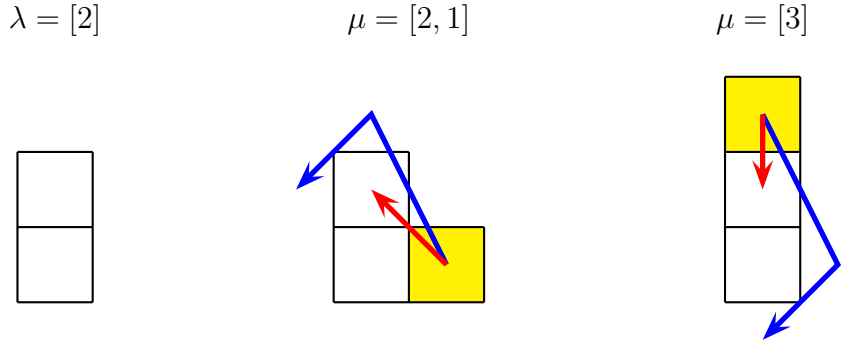
\begin{array}{ccccc}

\lambda = [2]  &  \hspace{.5cm}  &  \mu = [2,1]  &  \hspace{.5cm}  &  \mu = [3] \\

\begin{pspicture*}(0,-0.5)(3,3.5)
\psline(1,2)(1,0) \psline(2,2)(2,0) \psline(1,2)(2,2)
\psline(1,1)(2,1) \psline(1,0)(2,0)
\end{pspicture*}

&&

\begin{pspicture*}(0,-0.5)(4,3.5)
\psframe*[linecolor=yellow](2,1)(3,0) \psline(1,2)(1,0)
\psline(2,2)(2,0) \psline(3,1)(3,0) \psline(1,2)(2,2)
\psline(1,1)(3,1) \psline(1,0)(3,0)
\psline[linecolor=blue,linewidth=2pt,arrowsize=8pt]{->}(2.5,0.5)(1.5,2.5)(0.5,1.5)
\psline[linecolor=red,linewidth=2pt,arrowsize=8pt]{->}(2.5,0.5)(1.5,1.5)
\end{pspicture*}

&&

\begin{pspicture*}(0,-0.5)(3,3.5)
\psframe*[linecolor=yellow](1,3)(2,2) \psline(1,3)(1,0)
\psline(2,3)(2,0) \psline(1,3)(2,3) \psline(1,2)(2,2)
\psline(1,1)(2,1) \psline(1,0)(2,0)
\psline[linecolor=blue,linewidth=2pt,arrowsize=8pt]{->}(1.5,2.5)(2.5,0.5)(1.5,-0.5)
\psline[linecolor=red,linewidth=2pt,arrowsize=8pt]{->}(1.5,2.5)(1.5,1.5)
\end{pspicture*}

\end{array}
$$

\lpara

We deduce a formula for $q_{-1}$. Given a partition $\mu$, let $\mu[-1]$
denote the set
of partitions $\lambda$ with $\lambda \subset \mu$ and $|\lambda| = |\mu| - 1$.

\begin{prop}
\label{pro:q-1}
We have the following formula:
$$
q_{-1}(\fix(\mu)) = \sum_{\lambda \in \mu[-1]}
\frac{\Cok(\lambda,\mu)}{\Ker(\lambda,\mu)} \cdot
\frac{\Tan(\lambda)}{\Tan(\mu)}\ \fix(\lambda).
$$
\end{prop}
\begin{proo}
Let $\lambda$ resp. $\mu$ be a partition of weight $n$ resp. $n+1$.
Since $q_{-1}$ is the adjoint of $q_1$, we have the relation
$$
\scal{q_1(\fix(\lambda)),\fix(\mu)}_{S^{[n+1]}} =
\scal{\fix(\lambda),q_{-1}(\fix(\mu))}_\sn\ .
$$
Therefore the proposition follows from Lemma \ref{lem:norme-fix} and
Proposition \ref{pro:q1}.
\end{proo}

\subsection{Class of the boundary and derivatives}

\label{sub:boundary}

We turn to the problem of determining the equivariant class of the
divisor $\Delta_2$ of non-reduced schemes. In the non equivariant
setting on $S^{[n]}$, recall the following formula of Lehn which
expresses the class $[\Delta_2]_{cla}$
of $\Delta_2$ in terms of the classical Chern
class $c^{cla}$ of the tautological
bundle: $[\Delta_2]_{cla}=-2 c_1^{cla}({\cal O} ^{[n]})$. We prove an
equivariant analog in the equivariant Chow ring: $[\Delta_2]=-2
c_1({\cal O}^{[n]})$, where the Chern class considered is the
equivariant Chern class. Our method involves equivariant techniques
and does not rely on Lehn's ideas.

Denote $\partial \in A_T^*(S^{[n]})$ the class of $\Delta_2$, and let
$p:\Spec\ k \to S^{[n]}$ be a $T$-fixed point. We'd like to compute
$p^* \partial$. To this end, assume that $p$ corresponds to the partition
$\lambda = (\lambda_1 , \ldots , \lambda_l )$ of weight $n$.
We let $l(\lambda)$ denote
the number of non-vanishing parts of $\lambda$, and $h(\lambda) = \lambda_1$.
Let $\lambda^\vee$ denote the partition dual to $\lambda$.

\begin{prop}
\label{pro:diviseur}
We have
$$p^* \partial = - (\Sigma_{j=1}^{h(\lambda)} \lambda_j^\vee(\lambda^\vee_j-1))\,
U - (\Sigma_{i=1}^{l(\lambda)} \lambda_i(\lambda_i-1))\, V.$$
\end{prop}
\begin{proo}
We treat first the case $n=2$. Then $S^{[2]}$ is the
blow-up of $S \times S$ along the diagonal. Assume, moreover, that
$\lambda=(2)$. Then $T_p S^{[2]}$ contains 4 eigenlines, of weight
$-U,-U+V,-V,-2V$. In this case $\Delta$ is smooth, and
the tangent space $T_p \Delta$ contains the three eigenlines
of weight $-U,-V$ and $-U+V$: in fact, the first two lines are obtained
by translating the double point $p$, and $-U+V$ is the weight of the
deformation obtained with schemes supported at the origin. We deduce
that $p^* \partial = -2V$.

{\fl We} now consider the general case. Let $l = l(\lambda)$ and $h =
h(\lambda)$.
Given $\bx = (x_1,\ldots,x_l)$ and $\by = (y_1,\ldots,y_h)$
tuples of elements in $k$, we let $I_{( \bx , \by )}$
denote the ideal generated by the $l+1$ polynomials
$\prod_{i=1}^{m-1} (X-x_i) \cdot \prod_{j=1}^{\lambda_m} (Y-y_j)$,
where $m$ varies from $1$ to $l+1$. When all the $x_i$ and all the
$y_j$ are distinct, $k[X,Y]/I_{( \bx , \by )}$ is reduced and the
corresponding set of points is the set of $(x_i,y_j)$ where $i \leq
l$ and $j \leq \lambda_i$. Thus $I_{( \bx , \by )}$ has length $n$.
On the other hand, when $\bx = (0,\ldots,0)$ and $\by =
(0,\ldots,0)$, the ideal $I_{( \bx , \by )}$ is monomial and
generated by the elements the $X^mY^{\lambda_m}$, and thus also has
length $n$. Since the length of this family of ideals is
upper-semicontinuous, it follows that it is constant, and this
family is flat.

{\fl In this way}, we obtain a $T$-equivariant morphism
$\varphi : k^{l+h} \to S^{[n]}$
with respect to the natural action on $k^{l+h}$. 
We now compute $\varphi^* \partial$. If
$\{i_1,i_2\} \subset \{1,\ldots,l\}$ is a subset with two elements,
where we assume $i_1<i_2$, we denote by $\Delta_{\{i_1,i_2\}}
\subset k^{l+h}$ the class of the variety of tuples $(\bx,\by)$
with $x_{i_1} = x_{i_2}$ and $\partial_{\{i_1,i_2\}}$ its class in
the equivariant Chow ring of $k^{l+h}$. Let $z:\Spec k \to
k^{h+l}$ be the origin of $k^{h+l}$; since
$\partial_{\{i_1,i_2\}}$ is defined by one equation of weight $U$,
it follows that $z^*
\partial_{\{i_1,i_2\}} = -U$. Similarly,
if $j_1<j_2$, let $\Delta_{\{j_1,j_2\}}$
be the divisor defined by $y_{j_1} = y_{j_2}$, and let
$\partial_{\{j_1,j_2\}}$ denote its class. We have $z^*
\partial_{\{i_1,i_2\}} = -V$.

{\fl We} claim that
$$
\varphi^* \partial =
\Sigma_{\{i_1,i_2\} \subset \{1,\ldots,l\}} 2\lambda_{i_2}
\partial_{\{i_1,i_2\}} +
\Sigma_{\{j_1,j_2\} \subset \{1,\ldots,h\}} 2\lambda_{j_2}^\vee
\partial_{\{j_1,j_2\}}\ .
$$
Clearly, we have an equality of sets
$$\varphi^{-1}(\Delta)
= \bigcup_{\{i_1,i_2\} \subset \{1,\ldots,l\}} \Delta_{\{i_1,i_2\}}
\cup
\bigcup_{\{j_1,j_2\} \subset \{1,\ldots,h\}} \Delta_{\{j_1,j_2\}}\ ,
$$
and we claim that the multiplicity of $\Delta_{\{i_1,i_2\}}$ is
$2\lambda_{i_2}$. To see why, let $(\bx,\by)$ be a generic point in
$\Delta_{\{i_1,i_2\}}$: we have $x_{i_1} = x_{i_2}$ but no other equality
among the $x_i$'s and the $y_j$'s. Thus the scheme represented by
$\varphi(\bx,\by)$ is a union of $\lambda_{i_2}$ double points and
$n-2\lambda_{i_2}$ other distinct points.
Near the point
$\varphi(\bx,\by)$, $S^{[n]}$ is isomorphic to
$(S^{[2]})^{\lambda_{i_2}} \times S^{n-2\lambda_{i_2}}$.
Thus the multiplicity of our
component may be deduced from the case of $S^{[2]}$: in this case the
multiplicity was 2 in view of the computation we made at the beginning of
the proof. Thus the multiplicity is $2\lambda_{i_2}$ as claimed.

{\fl Since} $z^* \varphi^* \partial = p^* \partial$, it remains only to
show that
 $2 \Sigma_{\{i_1,i_2\} \subset \{1,\ldots,l\}} \lambda_{i_2}
= \Sigma_{j=1}^{h(\lambda)} \lambda_j^\vee(\lambda_j^\vee-1)$.
The first sum is equal to
$\Sigma_{1 \leq i_1<i_2 \leq l,1\leq j\leq \lambda_{i_2}}\ 2$. In this sum, when
$j=j_0$ is fixed, $i_2$ is such that $\lambda_{i_2} \geq j_0$,
which forces $i_2 \leq \lambda_{j_0}^\vee$. Thus
$\Sigma_{1 \leq i_1<i_2 \leq l,j\leq \lambda_{i_2},j=j_0}\ 2 =
\lambda_{j_0}^\vee(\lambda_{j_0}^\vee-1)$. Our proof is now complete.
\end{proo}

\begin{coro}
  \label{coro:big-diagonal}
  The equivariant class of $\Delta_2$ in $A_T^*(S^{[n]})$ is
  $\partial=-2c_1({\cal O}_X^{[n]})$.
\end{coro}
\begin{proo}
  By Proposition \ref{pro:diviseur}, the two classes have the same
  restriction on
  the $T$-fixed points of $S^{[n]}$  and the restriction morphism is injective.
\end{proo}

\lpara

If $f:A \to A$ is any operator, we now give a formula for the commutator
$[\partial,f]$. To express this formula, let us introduce the following
notation:

\begin{nota}
\label{nota:delta}
If $f:A \rightarrow A$ is an endomorphism, define
$\Delta_{f,\lambda,\mu} \in K$ for $\lambda,\mu$ partitions by the formula
\begin{displaymath}
f(\fix(\lambda))=
\sum _{\mu} \Delta_{f,\lambda,\mu}\fix(\mu).
\end{displaymath}
\end{nota}

For $c=(a,b) \in \N^2$, let $w(c)=aU+bV$ be the weight of the corresponding
monomial. Corollary \ref{coro:big-diagonal} immediately implies:

\begin{coro}
\label{coro:derivee}
Let $f:A \to A$ be any operator and let $\lambda \subset \mu$ be two partitions.
We have
$$
\Delta_{[\partial,f],\lambda,\mu} = -2 \Delta_{f,\lambda,\mu}\,
\sum_{c \in \mu \setminus \lambda} w(c)\ .
$$
\end{coro}

\subsection{Computation of the operator $q_i$ for all $i$}

In the previous sections, we computed $q_1,q_{-1}$ and $\partial$ on
the basis $\fix(\lambda)$. We introduce an auxiliary operator $\rho$
and give formulas for higher $q_i$'s in terms of $q_1,q_{-1}$ and
$\rho$. This yields an inductive procedure to compute $q_i$ on the
basis $\fix(\lambda)$.

\begin{defi}
\label{def:rho}
  Let $R^n \subset S^{[n,n+1]}$ be the closure of the
  set of pairs of schemes $(z_n,z_{n+1})$ with $z_n$ reduced,
  $z_n\subset z_{n+1}$ and $z_n = {(z_{n+1})}_{red}$.
\end{defi} 

Let  $\rho:\oplus_{n}\ A_T^*(S^{[n]}) \rightarrow \oplus_{n}\
A_T^{*+1}(S^{[n+1]}) $
be the morphism associated with the
correspondence $\amalg_n\ [R^n]$. It has conformal and cohomological
degree 1.

The following theorem gives a complete computation of the operators
$q_i$.
\begin{theo}
\label{theo:qi}
  We have
  $$
  \begin{array}{rlll}
  (i-1)q_i & = & \rho q_{i-1} - q_{i-1} \rho & \mbox{for } i>1 \\
  (i+1)q_i & = & \rho^\vee q_{i+1} - q_{i+1}\rho^\vee & \mbox{for } i<-1 \\
  2\rho & = & \partial q_1 -q_1 \partial \\
  2\rho^\vee & = & q_{-1} \partial - \partial q_{-1}
  \end{array}
  $$
\end{theo}
\begin{proo}
The non equivariant version of the first statement is proved
in \cite[Theorem 3.5]{lehn}.
Our formula can be proved geometrically as follows.
Let $\pi_1,\pi_2,\pi_3$ be the projections of
$S^{[n]} \times S^{[n+i-1]} \times S^{[n+i]}$ on each factor and, for
$a,b \in \{1,2,3\}$, let $\pi_{ab}$ be the projection on two factors.
To compute the composition $\rho q_{i-1}$, we have to understand the intersection
$\pi_{12}^{-1}(Q_{i-1}^n) \cap \pi_{23}^{-1}(R^{n+i-1})$. There are two
irreducible components in this
intersection. One, say $E_1$, is the closure of the set of triples of the form
$(z_n,z_n \amalg w_{i-1},z_{n+1} \amalg w_{i-1})$, where $z_n$ is a reduced subscheme of
length $n$, $w_{i-1}$ is a punctual subscheme of length $i-1$ with support not belonging
to $z_n$, and $z_{n+1}$ is a subscheme of length $n+1$ containing $z_n$ and having
the same support as $z_n$.

Another component denoted $E_2$
is the closure of the set of 
triples of the form 
$(z_n,z_n \amalg w_{i-1},z_n \amalg w_{i})$, 
where $z_n$ is again
a reduced subscheme of length $n$ and
$w_{i-1}$ resp. $w_{i}$ are
punctual subschemes
of length $i-1$
resp. $i$ with common support not belonging to $z_n$. The component
$E_2$ has multiplicity $i-1$ and $\pi_{13}(E_2) = Q_{i}^n$.
We claim that these are all the components of the intersection
$\pi_{12}^{-1}(Q_i^n) \cap \pi_{23}^{-1}(R^{n+i-1})$. This can be seen using arguments
similar to the detailed proof of Proposition \ref{pro:q1qi}; details
will be skipped here.

Consider now the composition $q_{i-1} \rho$ and the product
$S^{[n]} \times S^{[n+1]} \times S^{[i]}$. 
The intersection $\pi_{12}^{-1}(R^n) \cap \pi_{23}^{-1}(Q_{i-1}^{n+1})$ has only one
component $E'_1$ which is the closure of the set of triples
$(z_n,z_{n+1},z_{n+1} \amalg w_{i-1})$, with the same notations as for the component
$E_1$. In the commutator $\rho q_{i-1} - q_{i-1} \rho$ the components
$E_1$ and $E'_1$ cancel each other, and we get the formula.

\lpara

The third statement
is proved by a similar argument. The correspondences
in $S^{[n]} \times S^{[n+1]}$ corresponding to both
compositions $\partial q_1$ and $q_1 \partial$ contain the closure of the set of pairs
$(z_n,z_n \amalg w_1)$ where $z_n$ is a non-reduced subscheme of length $n$, and these
cancel each other. The composition $\partial q_1$ moreover contains the closure
of the set of pairs
$(z_n,z_{n+1})$ with $z_n$ reduced and $\supp(z_{n+1}) = \supp(z_n)$, namely, the
correspondence $R^n$, with mutliplicity 2.

\lpara

The second and the fourth equalities
are obtained from the first and the third equalities using duality and the fact that
$\partial$ is self-dual.
\end{proo}

\noindent
Applying this theorem and Corollary \ref{coro:derivee}, we deduce the following
formula for the operator $\rho$:

\begin{coro}
  \begin{equation}
    \label{equa:rho}
    \rho(\fix(\lambda)) = - \sum_{\mu \in \lambda[1]} \frac
    {\Cok(\lambda,\mu)}{\Ker(\lambda,\mu)}\, w(\mu \setminus \lambda)\, \fix(\mu) \ 
  \end{equation}
\end{coro}

\begin{exem}
\label{ex:nakToFix}
  Applying the result of this subsection recursivly, we obtain the
  following base
  changes between the basis $\nak(\lambda)$ and $\fix(\lambda)$ in
  conformal degree 2 and 3:
  \begin{eqnarray*}
\nak(1,1)&=&2\fix(1,1)+2\fix(2)\\
\nak(2)&=&-2U \fix(1,1)-2V \fix(2) 
\end{eqnarray*}

\begin{eqnarray*}
\nak(3)&=&6V^2 \fix(3)+ 3UV \fix(2,1)+6U^2 \fix(1,1,1)\\
\nak(2,1)&=&-6V \fix(3)-2(U+V) \fix(2,1) -6U \fix(1,1,1)\\
\nak(1,1,1)&=&6 \fix(3)+ 6 \fix(2,1)+6 \fix(1,1,1)
\end{eqnarray*}
The denominators of the fractions of the intermediate computations
simplify and the final base change is
polynomial. This is because $\nak(\lambda)$
lies in the subring $\ch{\sn}$ of $\ck{\sn}$.  
\end{exem}

\subsection{Commutation relations}

\label{sub:commutation}

In this subsection, we compute the commutators between the different
$q_i$'s.

We note that it is not possible to keep the proof by
Nakajima. Indeed, the equivariant pushforward of a class under a non
proper contracting morphism is not zero and the vanishing arguments of
Nakajima are not valid in our context. This non vanishing feature is
crucial for us because this is precisely the contribution of such
contracting morphisms that will give the non commutativity
$[q_{-1},q_{1}]=\frac{1}{UV} Id$.

\subsubsection{Commutation with $q_1$}
\label{sub:commWithq1}

 Our first goal is to study the commutator $[q_1,q_i]$. This
will follow from a geometric argument studying directly the
correspondences.

Recall (Definition \ref{def:qi}) that we denoted by $Q^n_i \subset
\sn \times S^{[n+i]}$ Nakajima's correspondence. Consider the
product $S^{[n]} \times S^{[n+1]} \times S^{[n+i+1]}$ and for $a,b
\in \{n,n+1,n+i+1\}$ the projection $\pi_{a,b}:S^{[n]} \times
S^{[n+1]} \times S^{[n+i+1]} \to S^{[a]} \times S^{[b]}$. Let us
denote by $\Ia \subset S^{[n]} \times S^{[n+1]} \times S^{[n+i+1]}$
the intersection $\pi_{n,n+1}^{-1}(Q_1^n) \cap
\pi_{n+1,n+i+1}^{-1}(Q_i^{n+1})$.












Let us introduce some piece of notation:

\begin{nota}
\label{nota:supp}
Let $w,z\subset S$ be two subschemes. Assume that $w \subset z$ or $w \supset z$.
If $w \subset z$ assume moreover that the support of ${\cal O}_z / {\cal O}_w$
is a point: in this case we denote by $\supp(w \neq z)$ this point.
If $w \supset z$ assume that the support of ${\cal O}_w / {\cal O}_z$
is a point: we denote by $\supp(w \neq z)$ this point.

Moreover, given a subscheme $z$ and a point $x$, we denote by $w_x$ the largest punctual
subscheme of $w$ whose support is $x$.

We denote by $l(w)$ the length of $w$.
\end{nota}

Let $E_1=\{(w,z,t)\in \Ia, z\cap t$ reduced, $supp(w\neq z)\neq
supp(z \neq t)\}$ and denote by $E_2$ the set $\{(w,z,t) \in \Ia,
z\cap t$ reduced, $supp(w \neq z)= supp(z \neq t)\}$.

\begin{prop}
\label{pro:q1qi}
The intersection $\Ia = \pi_{n,n+1}^{-1}(Q_1^n)
\cap \pi_{n+1,n+i+1}^{-1}(Q_i^{n+1})$ is proper. If $i>0$, then $\Ia
= \overline {E_1}$ and $\Ia$ is reduced irreducible of dimension
$2n+i+3$. If $i<0$, then $\Ia = \overline {E_1} \cup \overline{E_2}$
a union of two reduced subschemes of dimension $2n+i+3$.
\end{prop}
\begin{proo}
By Proposition \ref{pro:q1-lisse},
$\pi_{n,n+1}^{-1}(Q_1^n)$ is smooth and thus locally a complete
intersection. Therefore each irreducible component of $\Ia$ has
codimension at most $4n+i+1$ in $S^{[n]} \times S^{[n+1]} \times
S^{[n+i+1]}$, and so has dimension at least $2n+i+3$.

If $i<0$, let $e=2$ and if $i>0$, let $e=1$. To prove that $\Ia$ has
exactly $e$ reduced  components and the other claims of the
proposition, it suffices to describe a set of subschemes
$E(p,q)\subset \Ia$ and $E(p)\subset \Ia$ with the following
conditions:
\begin{itemize}
\item $\Ia=\coprod_{p,q} E(p,q) \amalg \coprod_p E(p)$ realizes $\Ia$ as a disjoint
union.
\item Exactly $e$ elements among the subschemes $E(p,q)$ and $E(p)$ have the
  expected dimension $2n+i+3$.
\item These $e$ strata are reduced.
\item The other strata have dimension less than $2n+i+3$.
\end{itemize}
The components in the intersection will then be the closures of the
maximal strata. For $p\geq 0$, $q\geq 0$, $i\neq 0$, $q+i\geq 0$,
let $E(p,q)$ be the set
\begin{displaymath}
\{(w,z,t)\in \Ia, \supp(w \neq z)\neq \supp(z \neq t),
l(w_{supp(w \neq z)})=p, l(w_{supp(z \neq t)})=q\}.
\end{displaymath}
For $p\geq 0$, $i\neq 0$, $p+1+i\geq 0$, let
\begin{displaymath}
E(p):=\{(w,z,t) \in \Ia, \supp(w \neq z)=supp(z \neq t)=x,
l(w_x)=p\}.
\end{displaymath}

\noindent Let $(w,z,t)$ in $E(p,q)$. Let $x=supp(w\neq z)$ and
$y=supp(z \neq t)$. Let $w_1\subset w$ the largest subscheme whose
support does not contain $x$ nor $y$. Since $w=w_1\cup w_x \cup
w_y$, $z=w_1 \cup z_x \cup w_y$, $t=w_1 \cup z_x \cup t_y$, the
triple $(w,z,t)$ is characterized by the data $w_1,
(w_x,z_x),w_y,t_y$. Since $l(w_1)=n-p-q$, $w_1$ moves in dimension
$2n-2p-2q$. The pair $(w_x,z_x)$ with $w_x\subset z_x$,
$l(w_x)=l(z_x)-1=p$ moves in dimension $p+2$ by Proposition \ref{prop-inco}.
The scheme $w_y$ with
$l(w_y)=q$ moves in dimension $q$ if $q=0$ and $q+1$ if $q>0$.
Given $w_y$, the
scheme $t_y$ with $l(t_y)=q+i$, $t_y\supset w_y$ (case $i>0$), $t_y
\subset w_y$ (case $i<0$) moves in dimension $q+i$ if $q+i=0$,
$q+i+1$ if $q+i>0$ and $q=0$, at most $q+i-1$ if $q+i>0$ and $q > 0$.

Summing up, in any case, the dimension of $E(p,q)$ is at most
$2n+i+3$, and the equality $\dim E(p,q)=2n+i+3$ is realized only
when $i>0$, $p=0$, $q=0$, and when $i<0,p=0,q=-i$.

\vskip 2mm

Let $(w,z,t)$ in $E(p)$. Let $x=supp(w\neq z)$ . Let $w_1\subset w$
the largest subscheme whose support does not contain $x$. Since
$w=w_1\cup w_x$, $z=w_1 \cup z_x$, $t=w_1 \cup t_x$, the triple
$(w,z,t)$ is characterized by the data $w_1, (w_x,z_x),t_x$. Since
$l(w_1)=n-p$, $w_1$ moves in dimension $2n-2p$. The pair $w_x,z_x$
with $w_x\subset z_x$ and $l(w_x)=l(z_x)-1=p$ moves in dimension
$p+2$. The scheme $t_x$ with $l(t_x)=p+1+i$, $t_x\supset z_x$ (case
$i>0$), $t_x \subset z_x$ (case $i<0$) moves in dimension $p+1+i$ if
$p+1+i=0$, at most $p+i$ if $p+1+i>0$.

Summing up, in any case, the dimension of $E(p)$ is at most
$2n+i+3$, and the equality $\dim S(p,q)=2n+i+3$ is realized only
when $i<0,p+1+i=0$.

By construction,  a point $(w,z,t)$ in a stratum of maximal
dimension is such that $z\cap t$ is reduced. The result follows.
\end{proo}

We now consider the product $S^{[n]} \times S^{[n+i]} \times
S^{[n+i+1]}$ and the three projections $\pi_{n,n+i}$,
$\pi_{n,n+i+1}$, $\pi_{n+i,n+i+1}$ defined as above. We denote by
$\Ib$ the intersection $\pi_{n,n+i}^{-1}(Q_i^n) \cap
\pi_{n+i,n+i+1}^{-1}(Q_1^{n+i})$.

Let $E'=\{(w,z,t)\in \Ia, z\cap t$ reduced, $supp(w\neq z)\neq
supp(z \neq t)\}$.

\begin{prop}
\label{pro:qiq1} The intersection $\pi_{n,n+i}^{-1}(Q_i^n) \cap
\pi_{n+i,n+i+1}^{-1}(Q_1^{n+i})$ is proper. More precisely $\Ib =
\overline {E'}$ and $\Ia$ is reduced irreducible of dimension
$2n+i+3$.
\end{prop}
\begin{proo}
The proof is similar to the proof of Proposition \ref{pro:q1qi}. We introduce a
stratification of $\Ib$ in the form
$\Ib=\coprod_{p,q} E(q,p) \amalg \coprod_q E(q)$.
For $p\geq 0$, $q\geq 0$, $i\neq 0$, $q+i\geq 0$, let $E(p,q)$ be the set

\begin{displaymath}
\{(w,z,t)\in \Ib, \supp(w \neq z)\neq \supp(z \neq t),
l(w_{supp(w \neq z)})=q, l(w_{supp(z \neq t)})=p\}.
\end{displaymath}
The only stratum of expected dimension $2n+i+3$ is $E(0,0)$ when
$i>0$ and $E(-i,0)$ when $i<0$. The study of the strata $E(p,q)$ is
rigorously similar to the mentioned Proposition  \ref{pro:q1qi} and
we skip it. For $q\geq 0$, $i\neq 0$, $q+i\geq 0$, let

\begin{displaymath}
E(q):=\{(w,z,t) \in \Ib, \supp(w \neq z)=supp(z \neq t)=x,
l(w_x)=q\}.
\end{displaymath}

Let $(w,z,t)$ in $E(q)$. Let $x=supp(w\neq z)$ . Let $w_1\subset w$
the largest subscheme whose support does not contain $x$. Since
$w=w_1\cup w_x$, $z=w_1 \cup z_x$, $t=w_1 \cup t_x$, the triple
$(w,z,t)$ is characterized by the data $w_1,(z_x,t_x),w_x$. Since
$l(w_1)=n-q$, $w_1$ moves in dimension $2n-2q$. The pair $(z_x,t_x)$
with $z_x\subset t_x$, , $l(z_x)=l(t_x)-1=q+i$ moves in dimension
$q+i+2$ by Proposition \ref{prop-inco}.
Given $z_x$, the scheme $w_x$ with $l(w_x)=q$, $w_x\subset z_x$ (case
$i>0$), $w_x \supset z_x$ (case $i<0$) moves in dimension $q$ if
$q=0$, at most $q-1$ if $q>0$.

Summing up, in any case, the dimension of $E(q)$ is at most
$2n+i+2$. There is no stratum $E(q)$ of the expected dimension
$2n+i+3$.
\end{proo}

In the proof of the next proposition we will use the following easy lemma:

\begin{lemm}
\label{lemm:bott}
Let $\pi : S \to pt$ be the projection of $S$ to a point. Then
$\pi_* 1 = 1/UV$.
\end{lemm}
\begin{proo}
This is a direct application of Theorem \ref{theo:bott}.
\end{proo}

\begin{prop}
  \label{pro:q1-qi}
  We have $[q_{-1},q_{1}]=\frac{1}{UV}Id$. Moreover,
  for $i\neq -1$, we have $[q_{i},q_1]=0$.
\end{prop}
\begin{proo}
The composition $q_iq_1$ resp. $q_1q_i$ corresponds to
$(\pi_{n,n+i+1})_* ( \pi_{n,n+1}^* [Q_1^n] \cup \pi_{n+1,n+i+1}^* [Q_i^{n+1}] )$ 
resp. $(\pi_{n,n+i+1})_* (\pi_{n,n+i}^* [Q_i^n] \cup \pi_{n+i,n+i+1}^* [Q_1^{n+1}] )$.

By Propositions \ref{pro:q1qi} and \ref{pro:qiq1}, the intersection
$\pi_{n,n+1}^{-1} (Q_1^n) \cap \pi_{n+1,n+i+1}^{-1} (Q_i^{n+1})$
resp. $\pi_{n,n+i}^{-1} (Q_i^n) \cap \pi_{n+i,n+i+1}^{-1}
(Q_i^{n+i})$ is proper and therefore the cup product $\pi_{n,n+1}^*
[Q_1^n] \cup \pi_{n+1,n+i+1}^* [Q_i^{n+1}]$ resp. $\pi_{n,n+i}^*
[Q_1^n] \cup \pi_{n+i,n+i+1}^* [Q_i^{n+1}]$ is equal to the class of
$\Ia$ resp. $\Ib$.

Moreover, when $i>0$, these two propositions show that $\Ia$ and
$\Ib$ are birational to $S^{[n]} \times S \times S^{[i]}_{punc}$,
where the indice ${punc}$ refers to the punctual Hilbert scheme.  So
there is a commutative $T$-equivariant diagram:
$$
\xymatrix{ \Ia \ar@{-->}[rr] \ar[rd]_{\pi_{n,n+i+1}} & &
\Ib \ar[ld]^{\pi_{n,n+i+1}} \\
& S^{[n]} \times S^{[n+i+1]} }
$$
From this it follows that, when $i>0$, the correspondence
$$
(\pi_{n,n+i+1})_* ( \pi_{n,n+1}^* [Q_1^n] \cup \pi_{n+1,n+i+1}^*
[Q_i^{n+1}] ) - (\pi_{n,n+i+1})_* ( \pi_{n,n+i}^* [Q_i^n] \cup
\pi_{n+i,n+i+1}^* [Q_1^{n+1}] )
$$
defining the commutator morphism $[q_1,q_i]$ is zero.

When $i<0$, there is an extra component in $\Ia$, namely
$\overline{E(-i-1)}$ with the notations of Proposition
\ref{pro:q1qi}. Let $c=[\overline{E(-i-1)}]$ denote the class of this
component. The commutator $[q_i,q_1]$ is then defined by the
correspondence $(\pi_{n,n+i+1})_*(c)$.

If $i<-1$, the morphism
$\pi_{n,n+i+1}:\overline{E(-i-1)}\fd S^{[n]} \times S^{[n+i+1]}$
 is a proper morphism with fibers of positive
dimension. It follows that the correspondence $(\pi_{n,n+i+1})_*(c)$
defining the morphism is equal to zero.

If $i=-1$,  the morphism
$\pi_{n,n+i+1}:\overline{E(-i-1)} \fd S^{[n]} \times S^{[n+i+1]}$
is not proper any more. It is birational
to the morphism $\varphi:\Delta \x S \fd S^{[n]} \times S^{[n]}$ where
$\Delta \subset S^{[n]} \times S^{[n]}$ is the diagonal and
$\varphi(w,w,x)=(w,w)$.
 It follows that the correspondence
$(\pi_{n,n+i+1})_*(c)$ defining the morphism is equal to
$\frac{1}{UV} [\Delta]$ and the proposition follows.
\end{proo}

\begin{prop}
\label{pro:commutativite-ij} Let $i,j$ be positive integers. Then
$q_iq_j = q_jq_i$.
\end{prop}
\begin{proo}
By Theorem \ref{theo:qi}, we have $(i-1)q_i = \rho q_{i-1} - q_{i-1}
\rho$ and $jq_{j+1} = \rho q_j - q_j \rho$. From this it follows
that
$$
(i-1)[q_i,q_j] = [ \rho , [q_{i-1},q_j] ] - j [ q_{i-1},q_{j+1} ]\ \
.
$$
By Proposition \ref{pro:q1-qi} and induction on $i$, we may assume
that $[q_{i-1},q_j] = 0$ and $[q_{i-1},q_{j+1}]=0$. Thus the
proposition is proved.
\end{proo}

\subsubsection{Commuting $\rho$ and $\rho^\vee$}
\label{sec:comm-relat}

We now compute the commutator $[\rho,\rho^\vee]$. This is the
technical key point of the computation of the commutation relations 
involving higher $q_i$'s.

\begin{prop}
\label{pro:rho-rhod}
  We have:
  $$[\rho,\rho^\vee]=\bigoplus_{n\geq 0}2n\ Id_{A_K^*(S^{[n]})}.$$
\end{prop}

The heart of the proof is to get rid of an
excess intersection component. To this aim, we use some standard
intersection theory formulas to break up the initial intersection product into
several pieces. After this rewriting, some of the intersections that
show up are transverse and
easy to compute. The other pieces (responsible for the excess
intersections) are intersections with Cartier divisors. They can be
handled with Chern class formalism.

\begin{proo}
Let us first compute the correspondence $\rho \rho^\vee$ in the
equivariant Chow ring of $S^{[n]} \times S^{[n]}$.

On the product $S^{[n]} \times S^{[n-1]} \times S^{[n]}$
we denote by $\pi_i$ and $\pi_{ij}$ ($i,j \in \{1,2,3\}$) the natural
projections. Let $C = \pi_{12}^{-1}({(R^\vee)}^n)$ and
$D = \pi_{23}^{-1}(R^{n-1})$. Let $E_1$ resp. $E_2$ in
$S^{[n]} \times S^{[n-1]} \times S^{[n]}$ be the closure of the set of triples
$(z_n,z_{n-1},z'_n)$ with $z_{n-1}$ reduced, $z_n$ and $z'_n$ non
reduced, $z_{n-1}\subset z_n$, $z_{n-1}\subset z_n'$ and
$\supp(z_n \not = z_{n-1}) = \supp(z'_n \not = z_{n-1})$ resp.
$\supp(z_n \not = z_{n-1}) \not = \supp(z'_n \not = z_{n-1})$.
Let $F_i := \pi_{13}(E_i)$: a generic element in $F_1$
resp. $F_2$ is a couple $(z_n,z'_n)$ where $z_n$ and $z'_n$ have the
same support, both have exactly one double point and the double points
have the same resp. different support. The generic elements in
$E_1,E_2$ are depicted in the following picture:

$$
\hspace{13mm}
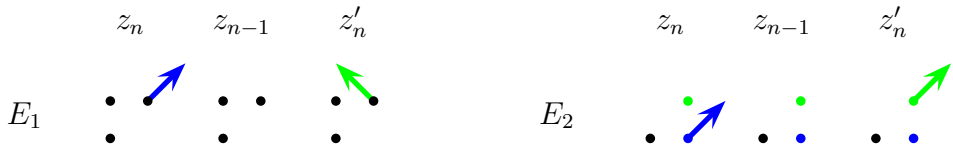
\begin{array}{cccccccc}

& && \hspace{-43mm} z_n \hspace{9mm} z_{n-1} \hspace{9mm} z'_n
& \hspace{-30mm} &&&
\hspace{-43mm} z_n \hspace{9mm} z_{n-1} \hspace{9mm} z'_n  \vspace{-17mm}\\

\comm{
& C &&
\begin{pspicture*}[shift=-.4](-.2,-.2)(8,3)
\psset{unit=5mm}
\psdots*(0,0)(0,1)(3,0)(3,1)(6,0)(6,1)
\psdots*[linecolor=blue](1,0)(4,0)(7,0)
\psdots*[linecolor=green](7,1)
\psline[linecolor=blue,linewidth=2pt,arrowsize=8pt]{->}(1,0)(2,1)
\end{pspicture*}

& \hspace{-30mm} &

D & \hspace{1mm} &

\begin{pspicture*}[shift=-.4](-.2,-.2)(8,3)
\psset{unit=5mm}
\psdots*(0,0)(0,1)(3,0)(3,1)(6,0)(6,1)
\psdots*[linecolor=blue](1,0)
\psdots*[linecolor=green](1,1)(4,1)(7,1)
\psline[linecolor=green,linewidth=2pt,arrowsize=8pt]{->}(7,1)(8,2)
\end{pspicture*}

\vspace{-17mm} \\
}

& E_1 &&
\begin{pspicture*}[shift=-.4](-.2,-.2)(8,3)
\psset{unit=5mm}
\psline[linecolor=blue,linewidth=2pt,arrowsize=8pt]{->}(1,1)(2,2)
\psline[linecolor=green,linewidth=2pt,arrowsize=8pt]{->}(7,1)(6,2)
\psdots*(0,0)(0,1)(1,1)(3,0)(3,1)(4,1)(6,0)(6,1)(7,1)
\end{pspicture*}

& \hspace{-30mm} &

E_2 & \hspace{1mm} &

\begin{pspicture*}[shift=-.4](-.2,-.2)(8,3)
\psset{unit=5mm}
\psdots*[linecolor=blue](1,0)(4,0)(7,0)
\psdots*[linecolor=green](1,1)(4,1)(7,1)
\psline[linecolor=blue,linewidth=2pt,arrowsize=8pt]{->}(1,0)(2,1)
\psline[linecolor=green,linewidth=2pt,arrowsize=8pt]{->}(7,1)(8,2)
\psdots*(0,0)(3,0)(6,0)
\end{pspicture*}

\\

\end{array}
$$

\begin{prop}
The intersection $C \cap D$ is generically transverse and
equal to the union $E_1 \cup E_2$.
\end{prop}
\begin{proo}
  The codimension of $C$ and $D$ in the product $S^{[n]} \times
  S^{[n-1]} \times S^{[n]}$ is $2n-1$. It follows that the components
  of $C\cap D$ have dimension at least $2n$. 
\comm{
\lt{Oui, c'est le fait que $U\cap V$ s'identifie a l'intersection de
  $U\times V$ avec la diagonale, ce qui donne la dimension attendue
  quand la variete ambiante est lisse car la diago est lisse
  egalement. C'est classique et je ne pense pas qu'il faille inclure
  l'argument dans le papier }
}

Moreover
$C\cap D\subset L$, where $L$ parametrizes the triples $(z_n,z_{n-1},z_n')$
  with $z_n\supset z_{n-1}, z'_n\supset z_{n-1},
  supp(z_n)=supp(z_{n-1}), supp(z_n')=supp(z_{n-1})$. The locus
  $L_k\subset L$ parametrizing the triples $(z_n,z_{n-1},z_n')$ with
  $z_{n-1}$ supported by $k$ points is such that $\overline{L_{n-1}}=E_1 \cup E_2$ has pure dimension 
$2n$. 
For $k<n-1$, $\dim L_k<2n$.
Thus the
  generic point of any component of $C\cap D$ is in $L_{n-1}$. The reverse
  inclusion $L_{n-1}\subset C\cap D$ is obvious, so that $C\cap
  D=\overline{L_{n-1}}=E_1\cup E_2$. The intersection is proper since 
  both $E_1$ and $E_2$ have dimension
  $2n$.  The intersection is transverse along $L_{n-1}$, thus
  generically transverse. 
\end{proo}

Since the
restrictions of $\pi_{13}$ to $E_1$ and $E_2$ are birational on their
image, it follows from the proposition that
\begin{equation}
\label{equa:rho-rhoDual}
\rho \rho^\vee = [F_1] + [F_2]\ .
\end{equation}

\vskip 6mm

Now we compute $\rho^\vee \rho$. We use similar notations for
$\pi_i,\pi_{ij}$ on $S^{[n]} \times S^{[n+1]} \times S^{[n]}$,
and moreover we denote by $\eta_1,\eta_2$ the two
projections from $\sn \times \sn$ to $\sn$. First of all we consider
the variety 
\begin{displaymath}
Q' = \pi_{12}^{-1}(Q_1^n) \cap
\pi_{23}^{-1}(Q_{-1}^{n+1}).
\end{displaymath}
We want to prove that $Q'$ admits
two irreducible components. 
\begin{lemm}\label{lemmeDim}
  Let $L_k\subset S^{[k-1]}\times S^{[k]} \times S^{[k-1]}$ be the locus
  parametrizing the triples $(z_{k-1},z_k,z_{k-1}')$ with
  $z_{k-1}\subset z_k$, $z_{k-1}'\subset z_k$ and $z_k$ supported at
  the origin. Then for $k\geq 2$, $\dim L_k\leq 2k-3$.  
\end{lemm}
\begin{proo}
  The pair $(z_{k-1},z_k)$ moves in dimension at most $k-1$. When
  $z_{k-1}$ and $z_k$ are fixed, $z_{k-1}'$ moves in dimension at most
  $k-2$. 
\end{proo}

An element in $Q'$ is a triple $(z_n,z_{n+1},z'_n)$ with $z_n\subset z_{n+1}$
and $z_{n+1}\supset z'_n$. Let $Q'_2\subset Q'$ the closed locus
defined by the condition $z_n=z_n'$ and $\underline{Q}'_1$
the open locus defined
by the condition $z_n\neq z_n'$. Let $Q'_1$ be the closure of
$\underline{Q}'_1$.  
\begin{prop}\label{QPrimeDefiniParIntersectionTransverse}
  The varieties $Q'_1$ and $Q'_2$ are irreducible of dimension $\dim
  Q'_1=\dim Q'_2=2n+2$. 
  The irreducible components of $Q'$ are
  $Q'_1$ and $Q'_2$. Moreover, the intersection $Q' = \pi_{12}^{-1}(Q_1^n) \cap
\pi_{23}^{-1}(Q_{-1}^{n+1})$ is generically transverse. 
\end{prop}
\begin{proo}
  The claims concerning $Q'_2$ are true
  since $Q'_2$ is isomorphic to $S^{[n,n+1]}$ by projection on the
  first two factors. 

As for $ \underline Q'_1$, let us denote by $p=z_{n+1}\setminus z_n$ and
$p'=z_{n+1}\setminus z_n'$ the natural residual points defined by 
a triple  $(z_n,z_{n+1},z'_n)\in \underline Q'_1$.

The irreducibility of $S_0^{[k,k+1]}$ implies that the locus
$L_0\subset \underline{Q}'_1$ with $p\neq p'$
is irreducible of dimension $2n+2$.  

For $m\geq 1$, we define the locus $L_m$ in $\underline{Q}'_1$
by the conditions $p=p'$
and $length(z_{n+1})_p=m$. It follows from Lemma 
\ref{lemmeDim} that it has dimension at most
$(2m-1)+2(n+1-m)=2n+1$ when $m\geq 2$. When $m=1$, $L_1=\emptyset$. 

Since the codimension of the intersection is bounded by the sum of the
codimensions, the components of $Q'$ have dimension at least $2n+2$. 

By construction $\underline Q'_1=\cup_{m\geq 2}L_m \cup L_0$ and
$Q'=\underline Q'_1 \cup Q'_2$. The dimensions computed
above show that the generic points of
$Q'$ coincide with the generic points of $L_0$ and 
$Q'_2$. Moreover,
for $m\geq 2$, $L_m\subset \overline{L_0}$ otherwise there would be in $Q'$ a
component of dimension less than $2n+2$.  

The transversality of the intersection $ \pi_{12}^{-1}(Q_1^n) \cap
\pi_{23}^{-1}(Q_{-1}^{n+1})$ is easily verified at the
generic points of $Q_1'$ and $Q'_2$. 
\end{proo}

We denote by
$C'_1$ resp. $C'_2$ the closures of the sets of triples
$(z_n,z_{n+1},z'_n)$ where $\supp(z_{n+1} \not = z_n) \in z_n$ and
$\supp(z_{n+1} \not = z_n) \not = \supp(z_{n+1} \not = z'_n)$ resp.
$\supp(z_{n+1} \not = z_n) = \supp(z_{n+1} \not = z'_n)$.
We
denote by $D'_1$ resp. $D'_2$ the closures of the sets of
triples $(z_n,z_{n+1},z'_n)$ where $\supp(z_{n+1} \not = z'_n) \in z'_n$
and $\supp(z_{n+1} \not = z_n) \not = \supp(z_{n+1} \not = z'_n)$ resp.
$\supp(z_{n+1} \not = z_n) = \supp(z_{n+1} \not = z'_n)$. 
The varieties $Q'_1,Q'_2,C'_1,C'_2,D'_1,D'_2$ (as well as the
following varieties $E'_1,E'_2,E'_3,E'_4$) are
depicted in the following array :

$$
\hspace{-10mm}
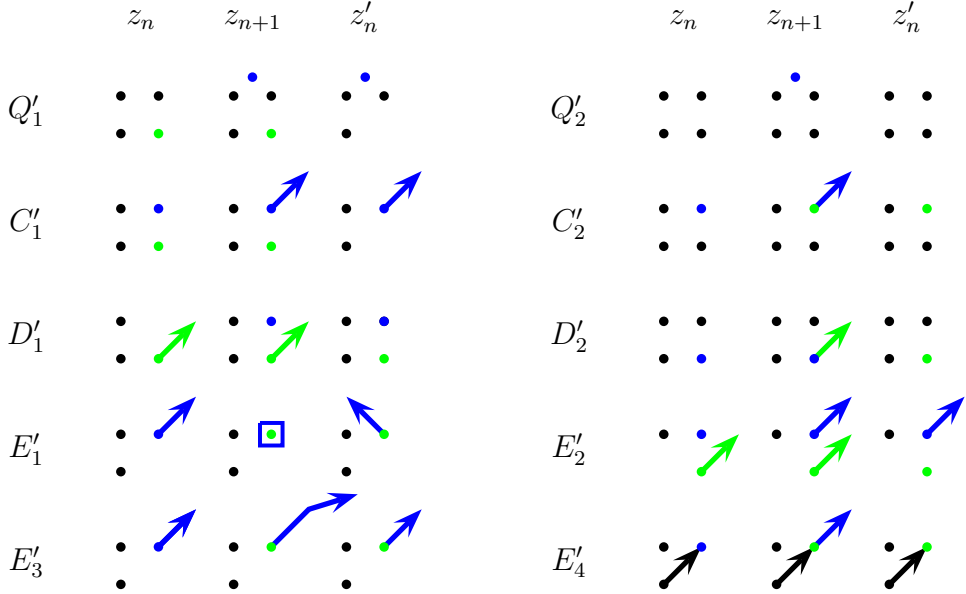
\begin{array}{cccccccc}

& && \hspace{-43mm} z_n \hspace{9mm} z_{n+1} \hspace{9mm} z'_n
& \hspace{-30mm} &&&
\hspace{-43mm} z_n \hspace{9mm} z_{n+1} \hspace{9mm} z'_n  \vspace{-17mm}\\

& Q'_1 &&
\begin{pspicture*}[shift=-.4](-.2,-.2)(8,3)
\psset{unit=5mm}
\psdots*(0,0)(0,1)(1,1)(3,0)(3,1)(4,1)(6,0)(6,1)(7,1)
\psdots*[linecolor=blue](3.5,1.5)(6.5,1.5)
\psdots*[linecolor=green](1,0)(4,0)
\end{pspicture*}

& \hspace{-30mm} &

Q'_2 & \hspace{1mm} &

\begin{pspicture*}[shift=-.4](-.2,-.2)(8,3)
\psset{unit=5mm}
\psdots*(0,0)(0,1)(1,0)(1,1)(3,0)(3,1)(4,0)(4,1)(6,0)(6,1)(7,0)(7,1)
\psdots*[linecolor=blue](3.5,1.5)
\end{pspicture*}

\vspace{-17mm} \\

\hspace{20mm} & C'_1 & \hspace{1mm} &

\begin{pspicture*}[shift=-.4](-.2,-.2)(8,3)
\psset{unit=5mm} \psdots*(0,0)(0,1)(3,0)(3,1)(6,0)(6,1)
\psdots*[linecolor=blue](1,1)(4,1)(7,1)
\psline[linecolor=blue,linewidth=2pt,arrowsize=8pt]{->}(4,1)(5,2)
\psline[linecolor=blue,linewidth=2pt,arrowsize=8pt]{->}(7,1)(8,2)
\psdots*[linecolor=green](1,0)(4,0)
\end{pspicture*}

& \hspace{-30mm} &

C'_2 & \hspace{1mm} &

\begin{pspicture*}[shift=-.4](-.2,-.2)(8,3)
\psset{unit=5mm}
\psdots*(0,0)(0,1)(1,0)(3,0)(3,1)(4,0)(6,0)(7,0)(6,1)
\psdot*[linecolor=blue](1,1)
\psline[linecolor=blue,linewidth=2pt,arrowsize=8pt]{->}(4,1)(5,2)
\psdots*[linecolor=green](4,1)(7,1)
\end{pspicture*}

\vspace{-17mm} \\

& D'_1 &&

\begin{pspicture*}[shift=-.4](-.2,-.2)(8,3)
\psset{unit=5mm} \psdots*(0,0)(0,1)(3,0)(3,1)(6,0)(6,1)(7,1)
\psline[linecolor=green,linewidth=2pt,arrowsize=8pt]{->}(1,0)(2,1)
\psline[linecolor=green,linewidth=2pt,arrowsize=8pt]{->}(4,0)(5,1)
\psdots[linecolor=green](1,0)(4,0)(7,0)
\psdots*[linecolor=blue](4,1)(7,1)
\end{pspicture*}

& \hspace{-30mm} & D'_2 &&

\begin{pspicture*}[shift=-.4](-.2,-.2)(8,3)
\psset{unit=5mm}
\psdots*(0,0)(0,1)(1,1)(3,0)(3,1)(4,1)(6,0)(6,1)(7,1)
\psline[linecolor=green,linewidth=2pt,arrowsize=8pt]{->}(4,0)(5,1)
\psdots*[linecolor=blue](1,0)(4,0)
\psdots*[linecolor=green](7,0)
\end{pspicture*}

\vspace{-17mm} \\

& E'_1 &&

\begin{pspicture*}[shift=-.4](-.2,-.2)(8,3)
\psset{unit=5mm} \psdots*(0,0)(0,1)(3,0)(3,1)(6,0)(6,1)
\psline[linecolor=blue,linewidth=2pt,arrowsize=8pt]{->}(1,1)(2,2)
\psline[linecolor=blue,linewidth=1.5pt](3.7,1.3)(3.7,0.7)(4.3,0.7)(4.3,1.3)(3.7,1.3)
\psline[linecolor=blue,linewidth=2pt,arrowsize=8pt]{->}(7,1)(6,2)
\psdot*[linecolor=blue](1,1) \psdots*[linecolor=green](4,1)(7,1)
\end{pspicture*}

& \hspace{-30mm} & E'_2 &&

\begin{pspicture*}[shift=-.4](-.2,-.2)(8,3)
\psset{unit=5mm}
\psdots*(0,1)(3,1)(6,1)
\psdots*[linecolor=blue](1,1)(4,1)(7,1)
\psdots*[linecolor=green](1,0)(4,0)(7,0)
\psline[linecolor=blue,linewidth=2pt,arrowsize=8pt]{->}(4,1)(5,2)
\psline[linecolor=blue,linewidth=2pt,arrowsize=8pt]{->}(7,1)(8,2)
\psline[linecolor=green,linewidth=2pt,arrowsize=8pt]{->}(1,0)(2,1)
\psline[linecolor=green,linewidth=2pt,arrowsize=8pt]{->}(4,0)(5,1)
\end{pspicture*}

\vspace{-17mm} \\

& E'_3 &&

\begin{pspicture*}[shift=-.4](-.2,-.2)(8,3)
\psset{unit=5mm} \psdots*(0,0)(0,1)(3,0)(3,1)(6,0)(6,1)
\psline[linecolor=blue,linewidth=2pt,arrowsize=8pt]{->}(1,1)(2,2)
\psline[linecolor=blue,linewidth=2pt,arrowsize=8pt]{->}(4,1)(5,2)(6.3,2.4)
\psline[linecolor=blue,linewidth=2pt,arrowsize=8pt]{->}(1,1)(2,2)
\psline[linecolor=blue,linewidth=2pt,arrowsize=8pt]{->}(7,1)(8,2)
\psdots*[linecolor=green](4,1)(7,1) \psdot*[linecolor=blue](1,1)
\end{pspicture*}

& \hspace{-30mm} & E'_4 &&

\begin{pspicture*}[shift=-.4](-.2,-.2)(8,3)
\psset{unit=5mm}
\psline[linecolor=blue,linewidth=2pt,arrowsize=8pt]{->}(4,1)(5,2)
\psline[linewidth=2pt,arrowsize=8pt]{->}(0,0)(1,1)
\psline[linewidth=2pt,arrowsize=8pt]{->}(3,0)(4,1)
\psline[linewidth=2pt,arrowsize=8pt]{->}(6,0)(7,1)
\psdots*(0,0)(0,1)(3,0)(3,1)(6,0)(6,1)
\psdots*[linecolor=blue](1,1)
\psdots*[linecolor=green](4,1)(7,1)
\end{pspicture*}

\end{array}
$$

Intersecting with Chern classes of line bundles, in particular with
Cartier divisors, commutes with the intersection product in the Chow ring
(\cite{fulton}, Example 8.1.6). It follows that if $X$ and $Y$ are
smooth in $Z$ smooth, if $\Delta_X\subset X$ and $\Delta_Y\subset Y$ are Cartier
divisors with restrictions $R_{X}$ and $R_{Y}$ on the (not necessarily smooth)
generically transverse intersection $X\cap Y$, then the intersection
$[\Delta_X]._Z[\Delta_Y]$ 
computed in the Chow ring $A^*(Z)$ is equal to $i_* (R_X._{X\cap Y}R_Y)$ where
$i:X\cap Y \fd  Z$ is the natural injection. Since we are working with
divisors, one can replace smoothness of $X$ and $Y$ with smooth in
codimension one. Moreover, there are equivariant analogs of these
statements. 

According to Proposition \ref{QPrimeDefiniParIntersectionTransverse}, 
one can apply 
the above with $Z=S^{[n]}\times S^{[n+1]} \times S^{[n]}$,
$X=\pi_{12}^{-1}(Q_1^n)$, $Y=\pi_{23}^{-1}(Q_{-1}^{n+1})$, $X\cap
Y=Q'=Q'_1\cup Q_2'$, $\Delta_X=\pi_{12}^{-1}(R^n)$,
$\Delta_Y=\pi_{23}^{-1}((R^\vee)^{n+1})$.
For the restrictions of the divisors, we use the notation 
$R_{Xi}=\Delta_X\cap Q_i$ and
$R_{Yi}=\Delta_Y\cap Q_i$. 
We obtain:
\begin{displaymath}
  \rho^\vee \rho =(\pi_{13})_*i_*( (R_{X1}+R_{X2}).(R_{Y1}+R_{Y2})) 
\end{displaymath}
where the intersection product takes place in $Q'$.

\begin{prop}
  $R_{X1}=[C'_1]$, $R_{X2}=2[C'_2]$, $R_{Y1}=[D'_1]$, $R_{Y2}=2[D'_2]$. 
\end{prop}
\begin{proo}
  This is clear set theorically. The multiplicities are computed in
  local coordinates at a generic point. 
\end{proo}

Since $R^n$ is a divisor on the smooth variety $S^{[n,n+1]}$,
$\Delta_X$ is a
Cartier divisor, and so are its restrictions $R_{Xi}$. Thus $[C'_1]$ and
$2[C'_2]$ are Cartier divisors, and similarly for $D'_1$ and $D'_2$.
Thus our task now is to compute the product of the divisors
$([C'_1]+2[C'_2])\cdot([D'_1]+2[D'_2])$ in $A_T^*Q'$. 
Note however that $C'_2=D'_2$, so that the
corresponding intersection is certainly not proper. In fact we
compute $2[C'_2]\cdot([D'_1]+2[D'_2])$ by another method. We know
that $\rho^\vee = \frac 12 (q_{-1}\delta - \delta q_{-1})$. So in
$A_T^*(Q')$ we have $2([D'_1] + 2[D'_2]) = \pi_2^* \partial_{n+1} -
\pi_3^* \partial_n$. The pushforward
$\pi_{13,*}([C'_2] \cup \pi_3^* \partial_n)$
can be computed thanks to the projection formula: this is
$\pi_{13,*} [C'_2] \cup \eta_2^* \partial_n$. But since $\pi_{13}$
is proper and contractant when restricted to $C'_2$,
$\pi_{13,*}\, [C'_2] = 0$.

To compute $\pi_{13,*}([C'_2] \cup \pi_2^* \partial_{n+1})$ we
observe that the general fibers of $\pi_{13}$ over $\pi_{13}(C'_2)$ are isomorphic $n$
copies of $\p^1$ and $\pi_2^* \partial_{n+1}$ restricts to a line
bundle isomorphic to ${\cal O}(-2)$ on each $\p^1$ (in fact the
class of the diagonal is $-2c_1^T({\cal O}^{[2]})$ if ${\cal O}^{[2]}$
denotes the tautological
bundle). Thus we get $\pi_{13,*}([C'_2] \cup \pi_3^* \partial_n) =
-2n Id_{\sn}$.

\medskip

To compute the other products we consider geometric intersections.
Let $E'_1,E'_2,E'_3,E'_4$ be the closures of some sets of triples
$(z_n,z_{n+1},z'_n)$. To define these triples we use the following
conventions: $p_i,p'_i$ will be punctual subschemes of length $i$ and
$w_j$ will be reduced subschemes of length $j$. Moreover, unless otherwise
stated, these subschemes will be generic (among punctual subschemes)
and their supports disjoint.

Let $p_2,p'_2$ share the same support,
and let $p_3$ be the 2-fat point having the same support as $p_2$
and $p'_2$.
A generic triple $(z_n,z_{n+1},z_n')$ in $E'_1$ is given
as follows: $z_n = w_{n-2} \amalg p_2\ ,\ z'_n = w_{n-2} \amalg p'_2$ and
$z_{n+1} = w_{n-2} \amalg p_3$.

Let $p_1\subset p_2, p'_1 \subset p'_2$, $w_{n-3}$ be generic.
A generic triple $(z_n,z_{n+1},z_n')$ in $E'_2$ is given
as follows: $z_n = w_{n-3} \amalg p_1 \amalg p'_2,
z_{n+1} = w_{n-3} \amalg p_2 \amalg p'_2,
z'_n = w_{n-3} \amalg p_2 \amalg p'_1$.

Let $p_2 \subset p_3$, $w_{n-2}$ be generic.
A generic triple
$(z_n,z_{n+1},z_n')$ in $E'_3$ is given as follows:
$z_n=z'_n=w_{n-2} \amalg p_2$ and $z_{n+1} = w_{n-2} \amalg p_3$.

Let $p_1 \subset p_2$, let $p'_2$ and $w_{n-3}$ be generic.
A generic triple $(z_n,z_{n+1},z_n')$ in $E'_4$ is given
as follows: $z_n = z'_n = w_{n-3} \amalg p_1 \amalg p'_2,
z_{n+1} = w_{n-3} \amalg p_2 \amalg p'_2$.

\begin{lemm}
\label{lem:inter-bords}
We have the set theoretic intersection $C'_1 \cap (D'_1 \cup D'_2)
= E'_1 \cup E'_2 \cup E'_3 \cup E'_4$.
\end{lemm}
\begin{proo}
Let $\cI$ denote an irreducible component in the intersection of the lemma.
We know that $\cI$ has dimension
at least $2n$. 
Let $\xi = (z_n,z_{n+1},z'_n)$ be a generic point
in one of these components.
Let $x = \supp(z_{n+1} \not = z_n)$ and $y = \supp(z_{n+1} \not = z'_n)$.
Let $p$ be the length of $z_n$ at $x$ and $q$ the length of $z'_n$ at $y$.
Since $\cI \subset C'_1$, $p \geq 1$. Since
$\cI \subset D'_1 \cup D'_2$, $q \geq 1$. Finally, since $\cI \subset C'_1$,
$z'_n$ is non reduced.

Assume first that $x \not = y$. Then the triple $\xi$ is defined by the
inclusions ${(z_n)}_{|x} \subset {(z_{n+1})}_{|x}$ and ${(z'_n)}_{|y} \subset {(z_{n+1})}_{|y}$
and the intersection
$z_n \cap (S \setminus \{x,y\})$ which has length $n-p-q-1$. So the dimension
of the set of such triples is $(p+2) + (q+2) + 2(n-p-q-1) = 2n+2-p-q$, so
that $p+q=2$. Therefore $p=1=q$, and so $\cI = E'_2$.

From now on, we assume that $x=y$, so $q=p$.
If $p=1$, since $z'_n$ is non reduced, we have $\cI = E'_4$.

Let us see that $p \leq 2$.
We denote by $f$ the dimension of the set of schemes of length
$p$ included in ${(z_{n+1})}_{|x}$. Since the support of such a subscheme of
${(z_{n+1})}_{|x}$ is $x$
we have $f \leq p-1$.
Let $r:\cI \dasharrow S_0^{[p,p+1]}$ which maps a triple $(z_n,z_{n+1},z'_n)$
to the pair $({(z_n)}_{|x},{(z_{n+1})}_{|x})$.
We denote by $d$ the dimension of
$r(\cI)$. By Proposition \ref{prop-inco} $\dim S_0^{[p,p+1]}=p+2$
so we have $d \leq p+2$.
Moreover $\dim \cI = 2n-2p+d+f\geq 2n$. Summing up, we have

\begin{equation}
\label{equ:dim-incidence}
d+f \geq  2p\ ,\ d \leq p+2\ ,\ f \leq p-1\ .
\end{equation}

If $f=p-1$, then any scheme of length $p$ supported at $x$ is included
in ${(z_{n+1})}_{|x}$, and this implies that $p=2$ and ${(z_{n+1})}_{|x}$ is a
$2$-fat point. In this case $\cI = E'_1$.

Let us assume that $f \leq p-2$. Equation (\ref{equ:dim-incidence}) implies
that $f=p-2$ and $d = p+2$, so $r(\cI) = S_0^{[n,n+1]}$.
If we assume that $p>2$, we get $f>0$ and therefore
${(z_{n+1})}_{|x}$ cannot be curvilinear, contradicting
$r(\cI) = S_0^{[n,n+1]}$.
Thus we have $p=2$ and $\cI = E'_3$.
\end{proo}

\begin{lemm}
\label{lem:c'-inter-d'}
There exist integers $a,b,c,d$ such that
$[C'_1] \cup [D'_1] = [E'_1] + [E'_2] + a[E'_3] + b[E'_4]$
and $[C'_1] \cup [D'_2] = c[E'_3] + d[E'_4]$.
\end{lemm}
\begin{proo}
Note that a generic point in $E'_1$ and $E'_2$ is a smooth point in
$Q'$ (this will for example be a consequence of our 
following parametrization
of $Q'$ near such a point).
Let us first compute the intersection number of $C'_1$ and $D'_1$ along
$E'_1$. A generic point in $C'_1$ resp. $D'_1,E'_1$ can be obtained by
disjoint union of $n-2$ distinct points and a generic point in the
same variety in the case $n=2$, thus it is enough to consider the case where
$n=2$. We consider the particular point $\xi = (z_2,z_3,z'_2)$ where
$z_2$ resp. $z_3,z'_2$ is the subscheme of the plane defined by the
equations $(X,Y^2)$ resp. $(X^2,XY,Y^2),(X^2,Y)$. Note that the projection
$S^{[2]} \times S^{[3]} \times S^{[2]} \to S^{[2]} \times S^{[2]}$ restricts
to an isomorphism on its image in a neighborhood of $\xi$ in $Q'$, since
for $\varepsilon = (y_2,y_3,y'_2)$ in such a neighborhood, $y_3$ is the
scheme-theoretic union of $y_2$ and $y'_2$. Thus $Q'$ is locally isomorphic
to the set of pairs $(y_2,y'_2)$ of subschemes of length $2$ which meet.
Note that for both $y_2$ and $y'_2$ there is a unique line containing it.
Moreover since all our intersection computations are invariant under
translations, we may assume that the intersection point of these two lines
is the origin.

We parameterize pairs of subschemes $(y_2,y'_2)$ near $(z_2,z'_2)$
such that these two lines meet at the origin by stating that $y_2$ resp. $y'_2$
corresponds to the ideal $(X+aY,Y^2+bY+c)$ resp.
$(X^2+dX+e,Y+fX)$. Then $Q'$ is defined by the fact that the origin
belongs to $y_2$ and $y'_2$, namely by the equations
$c=e=0$ (thus $Q'$ is locally an affine space).

Inside this variety, $C'_1$ resp. $D'_1$ is defined by the fact that
$y'_2$ resp. $y_2$ is non reduced. Thus it is defined by the equation
$d=0$ resp. $b=0$. We thus see that the intersection of $C'_1$ and $D'_1$
is transverse along a generic point in $E'_1$.

\medskip

Around a generic point in  $E'_2$
things are easier because the projection
$Q' \to S^{[n+1]},\xi \mapsto z_{n+1}$ is locally an isomorphism. 
Thus $Q'$ is locally
isomorphic to the product $S^{[n-3]} \times S^{[2]} \times S^{[2]}$, and
$(z_{n-3},z_2,z'_2)$ in this product belongs to
$C'_1$ resp. $D'_1$ if and only if $z_2$ resp. $z'_2$ is punctual.
So the intersection
$C'_1 \cap D'_1$ is transverse at such a point.

\medskip

Now the lemma follows from Lemma \ref{lem:inter-bords}.
\end{proo}

We have
$\pi_{13,*}\, [E'_3] = \pi_{13,*}\, [E'_4] =0$
since the restriction of $\pi_{13}$ to
$E'_3$ and $E'_4$ is proper and contractant. We have $\pi_{13,*}[E'_1]=[F_1]$
and $\pi_{13,*}[E'_2]=[F_2]$.
Therefore this gives $\rho^\vee \rho = [F_1] + [F_2] - 2n\, Id$. Since by
(\ref{equa:rho-rhoDual}), $\rho \rho^\vee = [F_1] + [F_2]$, 
the proposition is proved.
\end{proo}

\subsubsection{The commutator $[q_i,q_j]$}
\label{sec:commutator-q_i-q_j}

We can now compute the commutator $[q_i,q_j]$ for all $i,j$. 

\begin{lemm}
\label{lem:q-1-rho} We have $[q_{-1},\rho] = 0$.
\end{lemm}
\begin{proo}
First let us compute the correspondence $\rho q_{-1}$.
Consider the product $\sn \times S^{[n-1]} \times \sn$
and the natural projections
on this product.
Let $C := \pi_{12}^{-1}(Q_{-1}^n) \cap \pi_{23}^{-1}(R^{n-1})$. It is the
closure of the set of triples $(z_n,z_{n-1},z'_n)$ with $z_n$ reduced,
$z_{n-1} \subset z_n \cap z'_n$, and $z'_n$ having a point of length 2.
Let $F \subset \sn \times \sn$ be the closure of the set of pairs
$(z_n,z'_n)$ with $z_n$ reduced, ${(z'_n)}_{red} \subset z_n$ and $z'_n$
having a point of length 2 and simple points otherwise. Since the restriction of $\pi_{13}$ to
$C$ is birational, the morphism $\rho q_{-1}$ is given by the correspondence
$F$.

Now we compute the correspondence $q_{-1} \rho$.
The corresponding intersection has
been studied in the proof of Proposition \ref{pro:rho-rhod}.
With these notations
we have $\pi_{23}^{*}[Q_{-1}^{n+1}] \cup \pi_{12}^{*}[R^n] = [C'_1] + 2[C'_2]$.
Moreover the restriction of $\pi_{13}$ to $C'_1$ is birational with image $F$ and
the restriction of $\pi_{13}$ to $C'_2$ is proper contractant. Thus the morphism
$q_{-1} \rho$ is also given by the correspondence $F$, and the lemma is proved.
\end{proo}

Recall the convention that $q_0=0$.

\begin{prop}
\label{pro:qi-rho} Let $i$ be arbitrary. We have $[\rho,q_i] =|i|\,
q_{i+1}$.
\end{prop}
\begin{proo}
If $i \geq 0$ this is Theorem \ref{theo:qi}. If $i=-1$ this is Lemma
\ref{lem:q-1-rho}. Let us assume that $i=-j$ with $j \geq 2$. The
Jacobi identity reads:
$$
[[q_{-j+1},\rho^\vee],\rho] \ + \ [[\rho,q_{-j+1}],\rho^\vee] \ + \
[[\rho^\vee,\rho] , q_{-j+1}] = 0.
$$
Theorem \ref{theo:qi} yields
$[q_{-j+1},\rho^\vee] = (j-1) q_{-j}$. We may assume by induction
that $[\rho,q_{-j+1}] = (j-1) q_{-j+2}$. Finally by Proposition
\ref{pro:rho-rhod} we have $[[\rho^\vee,\rho] , q_{-j+1}] =
2(j-1)q_{-j+1}$. Thus we get:
$$
(j-1)[q_{-j},\rho] + (j-1)(j-2)q_{-j+1} + 2(j-1) q_{-j+1} = 0\ ,
$$
hence the proposition is proved.
\end{proo}

\begin{theo}
  \label{theo:qi-qj}
  Let $i$ and $j$ be any integers. We have
$$
[q_i,q_{-j}]\ =\ \left \{
\begin{array}{ll}
0 & \mbox{if } i \not = j
\vspace{2mm} \\
\frac{i(-1)^{i+1}}{UV}Id & \mbox{if } i=j
\end{array}
\right .
$$
\end{theo}
\begin{proo}
Since $q_{-i}$ is the adjoint of $q_i$ and $[q_i,q_{-j}] = 0$ if $i
\geq 0$ and $j \leq 0$ by Proposition \ref{pro:commutativite-ij}, we
may assume that $i,j \geq 1$. Moreover the proposition will be true
if $i=1$ or $j=1$ by Proposition \ref{pro:q1-qi}. Thus we assume
$i,j \geq 2$. Once again we apply Jacobi identity:
\begin{equation}
\label{equ:jacobi} [[\rho,q_{i-1}] , q_{-j}] \ + \
[[q_{-j},\rho],q_{i-1}] \ + \ [[q_{i-1},q_{-j}],\rho] \ = \ 0 \ .
\end{equation}
By induction we may assume that the commutator $[q_{i-1},q_{-j}]$ is
given by the proposition. Therefore it is either 0 or a scalar; in
both cases it will commute with $\rho$, so the last term vanishes.
By Proposition \ref{pro:qi-rho}, $[q_{-j},\rho] = -jq_{-j+1}$ and
$[\rho,q_{i-1}] = (i-1)q_i$.

Therefore equation (\ref{equ:jacobi}) reads $(i-1)[q_i,q_{-j}] =
j[q_{-j+1},q_{i-1}]$. If $i \not = j$, the second term vanishes by
induction and so $[q_i,q_{-j}] = 0$. If $i=j$ we get $[q_i,q_{-i}] =
i \frac{(-1)^{i+1}}{UV}Id $ as we wanted to prove.
\end{proo}

\section{Class of the small diagonal}
\label{sec:special-class}

Let $\Delta_i$ be the locus in $S^{[n]}$ where at least $i$ points
share the same support. In particular $\Delta_2$ is the
big diagonal, and $\Delta_n$ is the small diagonal. In Corollary
\ref{coro:big-diagonal},
we proved the equivariant formula for the big diagonal
$[\Delta_2]=-2c_1({\cal O}^{[n]})$, which is analogous to Lehn's formula
valid in the classical setting. In this section, we prove an
equivariant formula for the small diagonal.
\begin{theo}
\label{theo:classeSpeciale}
The $T$-equivariant class of
$\Delta_n$ is: $[\Delta_n]=(-1)^{n-1}\, n\, c_{n-1}({\cal O}^{[n]})$.
\end{theo}
The projection from the
equivariant Chow ring to the classical Chow ring gives obviously the
analogous formula in the classical setting.

\begin{rema}
\label{rema:lehn}
Given $u$ an equivariant line bundle over $S$,
let $c(u) \in \bigoplus A^*_T(\sn)$ denote
$([u^{[n]}])_n$, where $[\ \cdot \ ]$ denotes total equivariant Chern polynomial
and $u^{[n]}$ is the bundle over $\sn$ tautologically defined by $u$.
More generally we have the following formula:
$$
c(u) = \exp \left (
\sum_{m \geq 1} \frac{(-1)^{m-1}}{m} q_m \right ) \cdot \vacuum\ .
$$
\end{rema}
In this formula $\vacuum$ denotes the fundamental class in $S^{[0]}$.
\begin{proo}
This formula is Lehn's Theorem 4.6 \cite{lehn} and we explain why his proof
is valid in our equivariant context.
Lehn introduces the operator
$\mC(u) := c(u) \cdot q_1 \cdot c(u)^{-1}$ and shows
\cite[Theorem 4.2]{lehn} that
\begin{equation}
\label{equa:lehn}
\mC(u)=q_1(c(u)) + \rho \ .
\end{equation}
The proof of this
theorem relies on the exact sequence \cite[(11)]{lehn}
which is equivariant and his Lemma 3.9,
for which we proved an equivariant version (Corollary \ref{coro:derivee}).
Thus the relation (\ref{equa:lehn}) holds in the equivariant context.
Lehn's proof of \cite[Corollary 4.3]{lehn} is purely algebraic and therefore
we also have $c(u) = \exp(\mC(u)) \cdot \vacuum\ $. Finally, the proof of
\cite[Theorem 4.6]{lehn} uses this relation together with
the commutation relations
of $q_i$ and $\rho$, which we also proved in Theorem \ref{theo:qi}.
\end{proo}

We now give another proof of Theorem \ref{theo:classeSpeciale}, as a
straightforward consequence of an explicit expression of $q_n$
(Theorem \ref{theo:calculQn}) which we believe is interesting in itself.
The class $[\Delta_n]$ is equal to
$q_n \cdot \vacuum\ $.

Recall Notation \ref{nota:delta}.
If $\lambda \subset \N^2$ is a set of cardinal $n$ and
$M:\{1,\dots,n\}\fd \lambda$ is a bijection,
let $M^-:\{1,\dots,n-1\}\fd \lambda\setminus M(n)$ be the restriction of
$M$ and $M^+:\{1,\dots,n-1\}\fd \lambda\setminus M(1)$ the bijection
defined by $M^+(i)=M(i+1)$. Let $w:\lambda \fd \Q[U,V]$ be the
map sending $(a,b)$ to the linear form $aU+bV$
corresponding to the weight of the monomial $X^aY^b$ for the $T$-action.
Let
\begin{displaymath}
P_M=\frac{(-1)^{n-1}}{(n-1)!} \sum_{i=1}^{i=n}(-1)^{i-1} {\choose {i-1}{n-1}}
w(M(1))\cdots w(M(i-1))w(M(i+1))\cdots w(M(n))
\end{displaymath}
if $n>1$ and $P_M=1$ if $n=1$.

\begin{theo}
\label{theo:calculQn}
We have the relation
\begin{displaymath}
\Delta_{q_n,\lambda,\mu}=\sum_M
  P_M \, \prod_{i=0}^{n-1}\Delta_{q_1,\lambda_i,\lambda_{i+1}}\ \ ,
\end{displaymath}
where $M$ runs through the standard skew Young diagrams of shape
$\mu \setminus \lambda$, and  $\lambda_i$ is the partition defined by
$\lambda_i=\lambda \cup \{M(1),\ldots,M(i)\}$.
\end{theo}
Note that
$ \sum_M \, \prod_{i=0}^{n-1}\Delta_{q_1,\lambda_i,\lambda_{i+1}}
= \Delta_{q_1^n,\lambda,\mu}$.
If $\lambda$ is empty and $M$ is a tableau of shape $\mu$,
all the terms but the first
in the sum defining $P_M$ are zero, and
$P_M=\frac{(-1)^{n-1}}{(n-1)!} c_{n-1}({\cal O}^{[n]})_{|\fix(\mu)}$.
Since $q_1^n(\vacuum):=q_1\circ \cdots \circ
q_1(\vacuum)=n! \in A_T^*(S^{[n]})$, Theorem \ref{theo:classeSpeciale}
is a consequence of
Theorem \ref{theo:calculQn}.

\vskip 2mm

If $M:\{1,\dots,n\} \rightarrow \lambda$ is a standard skew young
diagram of shape $\lambda$, define $Q_M$ by $Q_M=1$ if $n=1$ and
recursively by the formula $Q_M=\frac{1}{n-1}
(-w(M(n))Q_{M^-}+w(M(1))Q_{M^+})$.
\begin{lemm}
  For every  standard skew Young diagram $M:\{1,\dots,n\} \rightarrow
  \lambda$,  $P_M=Q_M$.
\end{lemm}
\begin{proo}
  This is obvious if $n=1$ or $n=2$. To simplify the notation, we denote
  $w(M(k))$ by $m_k$. For $n$ general, we have
\begin{eqnarray*}
{(n-1)}Q_M&=&
-m_nQ_{M^-} + m_1Q_{M^+}\\
&=& -m_nP_{M^-} + m_1P_{M^+}\\
&=& -\frac{(-1)^{n-2}}{(n-2)!} (m_n \sum_{i=1}^{i={n-1}}(-1)^{i-1} \choose {i-1}
  {n-2}m_1\dots \hat m_i \dots m_{n-1}) \\
&\ &+m_1
  \sum_{i=2}^{i={n}}(-1)^{i}\choose {i-2 }
  {n-2}m_2 \dots \hat m_i \dots m_{n}) )\\
&=& \frac{(-1)^{n-1}}{(n-2)!} ( \sum_{i=2}^{i={n-1} }(-1)^{i-1}m_1\dots \hat m_i
\dots m_{n} \left ( \choose {i-1 }{n-2}+ \choose {i-2}{ n-2} \right ) \\
&\ & + m_2\dots
m_n+(-1)^{n-1}m_1\dots m_{n-1} )\\
&=& \frac{(-1)^{n-1}}{(n-2)!} \sum_{i=1}^{i={n}} (-1)^{i-1} \choose {i-1}
    {n-1}m_1\dots\hat m_i\dots m_{n}\\
&=& (n-1)P_M
\end{eqnarray*}
\end{proo}

\begin{lemm}
  \label{lem:coef-rho}
  Let $\lambda$ and $\mu$ be two Young diagrams of
  cardinal $n$ and $n+1$ with $\lambda \subset \mu$. Then
  $\Delta_{\rho,\lambda,\mu}=-w(\mu \setminus \lambda)
  \Delta_{q_1,\lambda,\mu}$.
\end{lemm}
\begin{proo}
  This is a direct consequence of the formula $2\rho=\partial q_1
  -q_1 \partial$ and the formula for $\partial$ given in Proposition
  \ref{pro:diviseur}.
\end{proo}

We now prove the formula for
$\Delta_{q_n,\lambda,\mu}$ from Theorem \ref{theo:calculQn}.
The formula is clearly true for $n=1$. Suppose that the formula for
$q_{n-1}$ is true. Since $(n-1)q_n=\rho q_{n-1}-q_{n-1}\rho$, we get:
\begin{eqnarray*}
  (n-1)\Delta_{q_n,\lambda,\mu}&=&\sum_{p_n\in Corners(\mu)}
  \Delta_{q_{n-1},\lambda,\mu\setminus{p_n}}\Delta_{\rho,\mu \setminus{p_n},\mu}
  \\ &\ &-\sum_{p_1\in
    OutsideCorners(\lambda)\cap \mu}\Delta_{\rho,\lambda,\lambda \cup {p_1}}
  \Delta_{q_{n-1},\lambda \cup \{p_1\},\mu}\\
  &\stackrel{\mbox{\footnotesize Lemma \ref{lem:coef-rho}}}{=}& \sum_{p_n\in Corners(\mu)}
  -\Delta_{q_{n-1},\lambda,\mu\setminus{p_n}}
  \Delta_{q_1,\mu \setminus{p_n},\mu}w(p_n)
  \\ &\ &+\sum_{p_1\in
    OutsideCorners(\lambda)\cap \mu}\Delta_{q_1,\lambda,\lambda \cup {p_1}}w(p_1)
  \Delta_{q_{n-1},\lambda \cup \{p_1\},\mu} \\
&\stackrel{\mbox{\footnotesize induction}}{=}& \sum_{M\mbox{ standard}}
  (\prod_{j=0}^{n-1}\Delta_{q_1,\lambda_j,\lambda_{j+1}})(-w(M(n))P_{M^-} + w(M(1))P_{M^+})\\
  &\stackrel{P_M=Q_M}{=}&(n-1)\sum_{M\mbox{ standard}}
  (\prod_{j=0}^{n-1}\Delta_{q_1,\lambda_j,\lambda_{j+1}})P_{M}
\end{eqnarray*}

\section{Base change formulas}

The goal of this section is to compute the base change formula from
$\es$ to $\nak$ and its inverse (recall Section \ref{section:bases} for the
bases $\es$ and $\nak$ of $A$). In particular, we prove that in the
classical setting, these two bases are equal up to a constant
(Theorem \ref{theo:nak2es}).

\subsection{Equivariant operators $q_{i,X}$}

The basis $\nak(\lambda)$ is defined using creation operators. The
basis $\es(\lambda)$ is defined via a Bialynicki-Birula
stratification. However, one can introduce operators $q_{i,X}$ such
that $\es(\lambda)$ is defined using creation operators too. The goal of
this section is to introduce the operators $q_{i,X}$ and to compute a
base change inductive formula  between $q_{i,X}$ and $q_i$
(Theorem \ref{theo:qix}).

The operator $q_{i,X}$ means ``adding $i$ points on a vertical line''. More
formally, $q_{i,X}:\ck{S^{[n]}} \to \ck{S^{[n+i]}}$
is defined by
the Fourier transform along the correspondence $Q_{i,X} \subset
S^{[n]} \times S^{[n+i]}$, where $Q_{i,X}$ is the closure of the set
of pairs $(z_n,z_n \amalg x_i)$ where $z_n \in S^{[n]},x_i \in S^{[i]}$,
$x_i$ is included in the vertical line $\Delta_{x_0}$ with equation $X=x_0$ for
some $x_0 \in k$, and $z_n$ and $x_i$ have disjoint support. We
denote by $\pi_{n}:Q_{i,X} \to S^{[n]}$ resp. $\pi_{n+i}:Q_{i,X} \to S^{[n+i]}$ the
natural projections.

First of all these operators allow the computation of the \ES cells:

\begin{prop}
\label{pro:op-ES}
Let $\lambda = (\lambda_1,\ldots,\lambda_l)$ be a partition. Then we have
$$
q_{\lambda_1,X} \circ \cdots \circ q_{\lambda_l,X} (\vacuum)
= \prod_i(\lambda^\vee_i- \lambda^\vee_{i+1})\, ! \cdot \es_\lambda\ .
$$
\end{prop}
\begin{proo}
To prove this result by induction on $l$, it is enough to show
that $q_{i,X}(\es_\lambda) = k\, \es_\mu$, where $\mu$ is the partition
obtained inserting one part equal
to $i$ in $\lambda$ and $k$ is the number of parts equal to $i$ in $\mu$.

To this end we apply the definition of $q_{i,X}$. For $n=|\lambda|$,
we have $q_{i,X}(\es_\lambda) = \pi_{n+i,*} \pi_n^*(ES_\lambda)$.
Recall that the Bialynicki-Birula cell decomposition by $ES_\lambda$ is
associated to the injection $k^* \to T, t \mapsto (t^{-1},t^{-d})$. 
Let $(z_n,z_{n+i})$ be a point belonging to
$\pi_n^{-1}(ES_\lambda)$, and assume that $z_n$ resp. $z_{n+i}$ belongs to the
open cell corresponding to the partition $\lambda'$ resp. $\mu'$. We claim
that $l(\mu') \leq l(\lambda)+1$.
In fact, since the whole construction is $k^*$-invariant, we also have
$(x_{\lambda'},x_{\mu'}) \in \pi_n^{-1}(ES_\lambda)$. On the other hand, for a generic
$(z_n,z_{n+i}) \in \pi_n^{-1}(ES_\lambda)$, the exists $x \in k$ such that
$(X-x) \cdot I(z_n) \subset I(z_{n+i})$ and thus we get
$X \cdot I(x_{\lambda'}) \subset I(x_{\mu'})$.
Therefore $l(\mu') \leq l(\lambda')+1$. Since $x_{\lambda'} \in ES_\lambda$,
$l(\lambda') \leq l(\lambda)$, thus $l(\mu') \leq l(\lambda)+1$.

Now, given such $\mu'$, we have $\dim (ES_{\mu'}) \leq n+l(\lambda)+i+1$.
In fact, the dimension of $ES_\lambda$ is equal to $n+l(\lambda)$.
Let $C$ be a component of $\pi_n^{-1}(ES_\lambda)$.
The dimension of $C$ is at least $n+l(\lambda)+i+1$. Thus, if the restriction
$C \to ES_{\mu'}$ is not dominant, it is contractant and it follows that
$\pi_{n+i,*} [C] = 0$. If it is dominant, then arguing on the generic points the
only possibility is that $\mu'=\mu$ and that $C$ is the component which is the
closure of the set of points $(z_n,z_{n+i})$ with $z_n$ generic in $ES_\lambda$
and $z_{n+i}$ obtained adding $i$ points on a vertical line to $z_n$. Let $C$
be this component.

The morphism $\pi_n:C \to ES_\lambda$ is submersive at a generic point in $C$,
so $C$ is a reduced component of $\pi_n^{-1}(ES_\lambda)$, thus
$\pi_{n+i,*} \pi_n^* [ES_\lambda] = \pi_{n+i,*} [ C ]$.

Moreover, given a generic element $z_{n+i} \in ES_\mu$, there are
$k$ vertical lines containing exactly $i$ points. Thus there
are $k$ couples $(z_n,z_{n+i})$ in the fiber $q^{-1}(z_{n+i})$:
$z_n$ is obtained from $z_{n+i}$ removing one of these lines.
Thus the restriction of
$\pi_{n+i}$ to $C$ has degree $k$ with image
$ES_\mu$, which proves the claim.
\end{proo}

Let $\Delta:=\Delta_0$ denote the vertical line with equation $X=0$.
Let $\SD{n}$ denote the subvariety of $S^{[n]}$ parameterizing
subschemes with support included in $\Delta$. If $\lambda$ is a
partition of weight $n$ and length $l$, let $\SDL{n}{\lambda}$
denote the closure in $\SD{n}$ of the variety of schemes
$z = z_1 \amalg \cdots \amalg z_l$, where $z_i$ has length $\lambda_i$ and is
supported on one point in $\Delta$.

\begin{prop}
  \label{pro:SD}
The varieties $\SDL{n}{\lambda}$ are the irreducible
components of $\SD{n}$, which is
therefore equidimensional of dimension $n$.
\end{prop}
\begin{proo}
Let $\lambda$ be a partition of weight $n$ and length $l$, and let
$i$ such that $1 \leq i \leq n$. The variety parameterizing schemes
of length $\lambda_i$ supported on one fixed point is irreducible by
\cite{briancon,ES}. Thus so is the variety parameterizing schemes of
length $\lambda_i$ supported on one point in $\Delta$. Thus each
$\SDL{n}{\lambda}$ is irreducible of dimension $n$. Since we have
$\SD{n} = \bigcup_\lambda \SDL{n}{\lambda}$, the proposition is
proved.
\end{proo}

\begin{theo}
\label{theo:qix}
We have the following formula:
$$
i\, q_{i,X} = (-1)^{i+1}\, q_i + U \cdot \sum_{j=1}^{i-1} (-1)^j \, q_j \circ q_{i-j,X}
$$
\end{theo}
With the help of
this theorem one can compute all the operators $q_{i,X}$ by induction on $i$.

\vskip 2mm

To prove the theorem we define auxiliary operators. For $i \geq 0$ and $j \geq 2$
let $q_{i,j,X}$ be the operator corresponding to ``adding $i$ points
on a same vertical line plus one punctual scheme of length $j$ whose
support is on this line''. Formally, $q_{i,j,X}$ is defined by an
incidence $Q_{i,j,X}$ in $S^{[n]} \times S^{[n+i+j]}$ where a
generic point in $Q_{i,j,X}$ is of the form $(z_n,z_n \amalg
\{x_1,\ldots,x_{i}\} \amalg t_j)$ where the $x_k$'s are distinct
points on a vertical line $\Delta$ not meeting $z_n$ and $t_j$ is a
length $j$ punctual scheme supported on $\Delta \setminus \{x_1,\ldots,x_{i} \}$.
Let us moreover use the convention that $q_{i,1,X} = (i+1) q_{i+1,X}$,
$q_{-1,j,X} = 0$, and $q_{0,X} = -1/U$.

\vskip 2mm

The theorem is a consequence of the following proposition
because this proposition implies that
the right hand side is equal to $q_{i-1,1,X} = i\, q_{i,X}$.

\begin{prop}
\label{pro:qijx}
For $i \geq 0$ and $j \geq 1$,
we have the following relation in $\ch{S^{[n]} \times S^{[n+i]}}$:
$$
-U \, q_j \circ q_{i,X} = q_{i,j,X} + q_{i-1,j+1,X}
$$
\end{prop}
\begin{proo}
If $i=0$ then the proposition is just formal thanks to the conventions
we made above. Let us assume that $i>0$.
Recall Definition \ref{def:qi}.
There is a morphism $s:Q_j \to S$ mapping the pair $(z_n,z_{n+j})$ to the support
of ${\cal O}_{z_{n+j}}/{\cal O}_{z_n}$. Let $\Delta$ be the line in
$S$ with equation $X=0$. Let $Q_j(-U) \subset Q_j$ denote
the divisor containing the set of pairs $(z_n,z_{n+j})$ in $Q_j$ with $s(z_n,z_{n+j}) \in \Delta$.
We have in $\ch{\sn \times S^{[n+j]}}$ the relation
$[Q_j(-U)] = s^* [\Delta] = -U \cdot [Q_j]$.
Similarly there is a morphism $x:Q_{i,X} \to \A^1$, mapping a pair
$(z_n,z_{n+i})$ to the common $X$-coordinate of the points in
$z_{n+i} \setminus z_n$.
Thus we use similar
notations and define $Q_{i,X}(-U) := x^{-1}(0) \subset Q_{i,X}$. In
the Chow ring, $[Q_{i,X}(-U)]=-U[Q_{i,X}]$.

Consider the product $\sn \times S^{[n+i]} \times S^{[n+i+j]}$ and
the projections $\pi_a,\pi_{ab}$. Let $I$ be the intersection
$\pi_{12}^{-1}(Q_{i,X}(-U)) \cap \pi_{23}^{-1}(Q_{j}(-U))$.
Every proper component of $I$ has dimension $2n+i+j$. 

Let $C$ be a component of $I$ which contributes to the composition 
$q_j \circ q_{i,X}$, ie. a component with $\pi_{13,*}[C]\neq 0$. 
Our first task is to prove that for a generic element
$(z_n,z_{n+i},z_{n+i+j})$
of $C$, the support of $z_n$ is disjoint
from $\Delta$.

If $z \subset S$ is a subscheme of dimension 0, we denote by $z_\Delta$ the
union of the components of $z$ supported on $\Delta$. Let $I_k$ be
the locally closed set of pairs $(z_n,z_{n+i+j})$ in $S^{[n]} \times
S^{[n+i+j]}$ such that the length of ${(z_n)}_\Delta$ is $k$, and $z_{n+i+j}
\supset z_n \cap (S \setminus \Delta)$,
the support
of ${\cal O}_{z_{n+i+j}} / {\cal O}_{z_n} \cap (S \setminus \Delta)$ 
is included in $\Delta$.
 Then $I_k$ is birational to
$S^{[k]}_\Delta \times S^{[n-k]} \times S^{[k+i+j]}_\Delta$, and
thus has dimension $2n+i+j$.

We denote by $k$ the integer such that for
a generic triple $(z_n,z_{n+i},z_{n+i+j})$ in $C$, the length
of ${(z_n)}_\Delta$ is $k$. Since, 
$\pi_{13}(C) \subset \overline I_k$, $\dim \pi_{13}(C)\leq 2n+i+j$.  
Moreover, if $k>0$, since $z_{n+i+j}$
has to contain $z_n$, $\pi_{13}(C)$ cannot contain $\overline I_k$, and
thus $\dim \pi_{13}(C) < 2n+i+j$. Since $\pi_{13}$ is proper we deduce
that $\pi_{13,*}[C] = 0$ in this case.

Let us now assume that $\dim \pi_{13}(C) = 2n+i+j$. We thus have $k=0$
and $\dim C=2n+i+j$. 
For a generic element $(z_n,z_{n+i},z_{n+i+j})$ in $C$,
${(z_{n+i+j})}_\Delta$ has length $i+j$, thus we have a
well-defined rational map $C \dasharrow S^{[i+j]}_\Delta$, with
$2n$-dimensional fibers. Let
$D$ be the closure of the image of this rational map.
Since $\dim D = i+j$, $D$ is a
component of $S^{[i+j]}_\Delta$; let us denote $\lambda$ the partition such
that $D = S^{[i+j]}_{\Delta,\lambda}$. By definition of
$I$, $\lambda$ must be dominated by the partition
$(j,1^i)$. It is clear that $D$ can contain $S^{[i+j]}_{\Delta,\mu}$
only if $\mu=(j,1^i)$ or $\mu=(j+1,1^{i-1})$. Therefore
$I$ has exactly two components which are not contracted
by $\pi_{13}$.

To describe these components let us consider some subschemes
$z_n,x_{i-1},x_i,p_j,p_{j+1}$ satisfying the following conditions.
The lengths of these subschemes are given by their indices.
The support of $z_n$ does not meet $\Delta$, whereas the other subschemes
have support included in $\Delta$. The subschemes $p_j,p_{j+1}$ are
punctual whereas $x_{i-1}$ and $x_{i}$ are reduced. Finally
$p_j \subset p_{j+1}$, $x_{i-1} \subset x_i$, and the support of
$p_{j+1}$ is not included in $x_i$. With these conditions
let $I_1$ resp. $I_2$ be the closure of the set
of triples $(z_n,z_{n+i},z_{n+i+j})$ where
$z_{n+i} = z_n \amalg x_i$ and
$z_{n+i+j} = z_n \amalg x_i \amalg p_j$ resp.
$z_{n+i+j} = z_n \amalg x_{i-1} \amalg p_{j+1}$.
The restriction of $\pi_{13}$ to $I_2$ is birational with image
$Q_{i-1,j+1,X}(-U)$. We have
$\pi_{13}(I_1)=Q_{i,j,X}(-U)$. If $j>1$ then the restriction
of $\pi_{13}$ to $I_1$ is birational whereas if $j=1$ it has degree
$i+1$. In view of our convention for $q_{i,1,X}$, this proves the proposition.
\end{proo}

\vskip 5mm





\subsection{Base change formulas}
\label{sec:base-change-formulas}

\begin{defi}
If $\lambda \in {\cal P}_n$ is a partition, we define
the operators $q_{\lambda}=\Pi_{i \in \lambda}q_i$,
$q_{\lambda,X}=\Pi_{i \in \lambda}q_{i,X}$, and the constant
$z_\lambda= \Pi \lambda_i \Pi
(\lambda_i^\vee-\lambda_{i+1}^\vee)! $. Let $j\in \lambda$.
With the notation with
multiplicity $\lambda=(1^{\alpha_1}, \dots, r^{\alpha_r})$, 
we denote by $\lambda \setminus j$ the
partition  $(1^{\alpha_1}, \dots,j^{\alpha_j-1},\dots,
r^{\alpha_r})$, with multiplicity one less for $j$. We let $t_\lambda=
\sum_{j\in \lambda} \frac{j}{(\alpha_j-1)!}  \frac{(l(\lambda)-1)!}
{\Pi_{i\neq j}(\alpha_i !)} $ and
$u_{\lambda}= \prod_i(\lambda^\vee_i - \lambda^\vee_{i+1})!$.
\end{defi} 

By definition of $q_{\lambda}$ and $q_{\lambda,X}$, the base
change formulas from $q_{\lambda}$ to $q_{\lambda,X}$ are determined
by the decomposition of $q_n$ in terms of the operators
$q_{\lambda,X}$ and similarly for the inverse base change. In
particular, the following theorem gives a full base change at the
level of operators. Since $\es(\lambda)=\frac{1}{u_\lambda}
q_{\lambda,X}(\vacuum)$ and  since $\nak(\lambda)=q_{\lambda}(\vacuum)$
the theorem applied to the vacuum also
yields the corresponding base
changes between $\es(\lambda)$ and $\nak(\lambda)$.

\comm{Formule verifiee a l'ordi}
\begin{theo}
\label{theo:chgt}
  $$ q_{i,X} = (-1)^{i+1}\, \sum_{|\lambda|=i } z_\lambda^{-1}\,
  U^{l(\lambda)-1}\, q_\lambda
  $$
 $$ q_{i}=(-1)^{i+1}\sum_{|\lambda|=i } t_\lambda \,
  U^{l(\lambda)-1}\, q_{\lambda,X}
  $$
\end{theo}
\begin{proo}
  By induction, the case $i=1$ being obvious.
  \begin{eqnarray*}
    i\, q_{i,X} & \stackrel{\mbox{\footnotesize Theorem \ref{theo:qix}}}{=} & (-1)^{i+1} q_i +
    U\, \sum_{j=1}^{i-1} (-1)^{j}\, q_j\circ q_{i-j,X} \\
  & \stackrel{\mbox{\footnotesize induction hypothesis}}{=}& (-1)^{i+1} q_i +
  \sum_{j=1}^{i-1} (-1)^{j}\, U\, (-1)^{i-j+1}\, \sum_{|\lambda|=i-j} \, z_{\lambda}^{-1}\,
  U^{l(\lambda)-1} \, q_j \circ q_{\lambda} \\
&=& (-1)^{i+1} q_i + (-1)^{i+1} \sum_{|\mu|=i, l(\mu) > 1} c_{\mu}\, q_{\mu}
\end{eqnarray*}
with
\begin{eqnarray*}
c_{\mu} & = & \sum_{j\in \mu} z_{\mu \setminus j}^{-1}\,
U^{l(\mu)-1} \\
\end{eqnarray*}
Since $\sum_{j\in \mu} \frac{z_\mu}{z_{\mu \setminus j}}=|\mu|$, we obtain
$c_{\mu}=|\mu|\, (-1)^{|\mu|+1}\, z_{\mu}^{-1}\, U^{l(\mu)-1}$, as required for
the induction.

The proof of the second formula is similar : the difficulty is to guess
the formula for $q_i$, then the induction is straightforward. 
Indeed, we start with the formula of theorem  \ref{theo:qix}
$(-1)^{i+1} q_i = -i q_{i,X}+
    U\, \sum_{j=1}^{i-1} (-1)^{j}\, q_j\circ q_{i-j,X}$,
and we replace $q_j$ on the right hand side by the induction formula. 
With the value
of $t_\lambda$ in the definition above and the
the formula for $q_i$, the induction follows.


\end{proo}

We now project the previous theorem from the equivariant Chow ring
to the classical Chow ring. All the constructions made so far
in the equivariant setting can be realized in the classical
setting. We denote by $q_{n}^{cla}$ and $q_{n,X}^{cla}$ the
corresponding operators on the classical Chow ring. Similarly, we
denote by $\nak^{cla}(\lambda)$ and $\es^{cla}(\lambda)$ the bases of
the classical Chow ring induced by these operators.

\begin{theo}\label{theo:nak2es}
  $$
  \begin{array}{rcl}
  q_n^{cla} & = & (-1)^{n+1}\, n\, q_{n,X}^{cla} \vspace{3mm}\\
  \nak^{cla}(\lambda) & = & (-1)^{|\lambda|+l(\lambda)}\, (\prod_{i\in
  \lambda}i )\, \es^{cla}(\lambda)
  \end{array}
  $$
\end{theo}
\begin{proo}
  In the classical setting, $U=0$ and the first formula for the
  operators is the
  projection of the corresponding formula in the equivariant
  setting. Applying the operators to the vacuum yields the second
  formula.
\end{proo}


\end{document}